\documentclass[preprint,12pt]{elsarticle}

 \usepackage{xcolor}
\usepackage{hyperref}

\usepackage[left=2cm,right=2cm,top=2cm,bottom=2cm]{geometry}

\usepackage{amssymb}

\usepackage{amsmath}
\DeclareMathAlphabet{\mathbfsf}{\encodingdefault}{\sfdefault}{bx}{n}
\numberwithin{equation}{section}

\usepackage{dsfont}
\usepackage{enumitem}

\usepackage[normalem]{ulem} 

\newlist{steps}{enumerate}{1}
\setlist[steps, 1]{label = Step \arabic*:}

\newtheorem{thm}{Theorem}
\newtheorem{corollary}[thm]{Corollary}
\newtheorem{lem}[thm]{Lemma}
\newtheorem{proposition}[thm]{Proposition}

\newdefinition{definition}{Definition}
\newdefinition{rmk}{Remark}
\newdefinition{exmp}{Example}
\newproof{pf}{Proof}
\newproof{pflemma1}{Proof of Lemma~\ref{lemma1}}
\newproof{pflemma2}{Proof of Lemma~\ref{lemma2}}

\newcommand  {\ZZ}{\mathbb{Z}}
\newcommand  {\N}{\mathbb{N}}

\newcommand {\corpsValeurs} {\mathbb{C}}
\newcommand {\corpsC} {\mathbb{C}}
\newcommand {\corpsR} {\mathbb{R}}

\newcommand{\vertiii}[1]{{\left\vert\kern-0.25ex\left\vert\kern-0.25ex\left\vert #1 
        \right\vert\kern-0.25ex\right\vert\kern-0.25ex\right\vert}}

\journal{---} 
\newcommand{\seb}[1]{#1}

\begin{document}

\begin{frontmatter}

\title{Exponential stability of linear periodic difference-delay equations} 
\date{August 31, 2024}

\author{L. Baratchart}
\author{S. Fueyo \footnote{Corresponding author.}}
\author{J.-B. Pomet}
\address{Inria, Université Côte d'Azur, Teams FACTAS and MCTAO,
  \\ 2004, route des Lucioles, 06902 Sophia Antipolis, France.}

\begin{abstract}
  This paper deals with the stability of linear periodic
  difference delay systems, where the value at time $t$ of a solution 
  is a linear combination with periodic coefficients
  of its values at finitely many delayed instants $t-\tau_1,\ldots,t-\tau_N$.
  We establish a necessary and sufficient condition for exponential stability
  of such systems when the coefficients have H\"older-continuous derivative,
  that generalizes  the one obtained  for difference delay systems with constant coefficients by
  Henry and Hale in the 1970s. This condition may be construed as
   analyticity, in a half plane,
  of the (operator valued) harmonic transfer function
    of an associated linear control system.
\end{abstract}


\end{frontmatter}

\section{Introduction}

In this paper we present a necessary and sufficient condition, stated in the frequency domain, for exponential stability of periodic
\emph{difference-delay systems}\footnote{
  There seems to be no general agreement on  terminology.
  The name ``difference-delay system'' is used in \cite{Ngoc2015,BFLP,SFueyo}, 
   but these are  called \emph{time-delay systems} in \cite{CoNg} and
  \emph{difference equations} in \cite{Henry1974,Chitour2016},
  while no specific name  is coined in \cite{Hale,CRUZ1970} where \eqref{systae}
   is written 
  $L\,y(.)=0$ and  $L$ is referred to as a \emph{difference operator}.
};
{\it i.e.}, linear dynamical systems of the form
\begin{equation}
\label{systae}
y(t)= \sum_{j=1}^ND_j(t)\,y(t-\tau_j), \qquad
t> s\,,
\end{equation}
where  $\tau_1<\cdots<\tau_N$ are
positive delays and $D_1(t),\ldots,D_N(t)$  complex $d\times d$ matrices,
depending  periodically on time $t$. 
This stability condition applies 
when the maps $t\mapsto D_j(t)$ are 
periodic and differentiable with Hölder continuous derivatives.
Here, periodicity is essential but the H\"older-smoothness assumption is technical.
Precise definitions of exponential stability are made in  Section~\ref{section_results};
for the time being, we simply describe
it as the property that every \emph{solution} 
  $y:[s-\tau_N,+\infty)\to\corpsValeurs^d$ of \eqref{systae} has a
  restriction to   $[t-\tau_N,t]$ that decays exponentially fast
with $t$ as the latter goes to $+\infty$.
\seb{In this framework, the main result is
  Theorem~\ref{theorem_hale_generalization} below, with a control
  theoretic counterpart in Theorem~\ref{theorem_hale_generalization2}
  formulated in terms of Harmonic Tranfer Function.}

\medskip

Dynamical systems like \eqref{systae} are
  natural generalizations of time-invariant difference-delay systems,
  arising  in various contexts of modeling and control.
Let us mention three:
\begin{enumerate}[label=\alph*)]
\item\label{motiv-pde}
  The study of certain \seb{(networks of)} \emph{1-D hyperbolic systems} reduces to the study of particular
  difference-delay systems: after integrating in terms of
  functions of one variable (backward and forward waves) using the
  method of characteristics, their evolution
  is governed by equations like \eqref{systae},
  see \seb{
  \cite{bastin2016stability,CoNg}, or the introduction of \cite{BFLP}.   
The results of the present paper can be applied to obtain necessary and sufficient stability criteria for
  certain 1-D hyperbolic PDE's (conservation laws) with linear periodic boundary \seb{conditions}, see the conference paper \cite{BARATCHART2022228}. \seb{In fact,
    the link with 1-D hyperbolic systems naturally arises {\it via}} the method of characteristics, allowing one to construe such systems as a linear periodic difference-delay systems of the form \eqref{systae}.}

\item\label{motiv-neutral}
  \emph{Neutral functional differential equations} are functional differential equations of the form
  \begin{equation}
    \label{eq:2}
    \textstyle\frac{\mathrm{d}}{\mathrm{d}t} \left(y(t)- \sum_{j=1}^N D_j(t) \,y(t-\tau_j)\right)
    = B_0(t)\,y(t)+\sum_{j=1}^NB_j(t)\,y(t-\tau_j),\qquad t\geq s,
  \end{equation}
  where the maps $B_i(.)$ and $D_i(.)$ are regular enough.
    Such equations model quite general linear phenomena involving
    delays, and are naturally related to  \eqref{systae}.
Indeed, the solution operator of System~\eqref{eq:2} turns out to be
a compact perturbation of the one governing System~\eqref{systae}, which
makes exponential stability properties of \eqref{eq:2} close
to those of \eqref{systae}: the spectra of their monodromy operators
differ by at most finitely many eigenvalues \cite[thm 7.3 Section 3.7]{Hale}.
References \cite{CRUZ1970,Henry1974} use this connection already
in the time-invariant case, and one can surmise that
characterizing the exponential stability of 
\eqref{systae} is likewise a substantial step towards analyzing the
one of \eqref{eq:2} in the periodic case;
see  discussion in Section \ref{sec:conj}.
In the language of electronic engineering, System~\eqref{systae} is called the \emph{high frequency limit system} of \eqref{eq:2} because it represents the limiting behaviour of the latter when $y$ oscillates arbitrary fast (so that the right hand side of \eqref{eq:2} goes to zero weakly).

\item\label{motiv-amplis}
  The stability of  \emph{active microwave circuits} (like amplifiers),
    has been
    an initial motivation of the authors to undertake the present study.
  Circuits can be regarded as nonlinear systems with infinite dimensional state space, due
to transmission lines that are modelled either by 1-D hyperbolic systems or
by delays. Their response to a periodic
signal
is typically a periodic solution, and linearizing along this trajectory yields
a linear periodic
dynamical system, governed by a functional differential equation
slightly more general 
than \eqref{eq:2},
whose operator solution is again a compact perturbation of a
  high-frequency limit system of the form \eqref{systae}.
The same
  remarks as in \ref{motiv-neutral} now apply, to the effect that
if the high frequency limit system is exponentially stable
(this is the property we shall characterize),
  then
the local stability of the periodic solution to the initial, nonlinear circuit depends only on whether
  the harmonic transfer function
of the linearized system (HTF, to be defined later) has unstable \emph{poles} \cite{SFueyo}.
\\
Such methods are already instrumental in engineering (without much mathematical
justification), based on frequency-wise simulation of the HTF obtained by so-called
\emph{harmonic-balance} techniques; see for example \cite{Suarez} for
a  discussion.
\seb{Here, it it important to stress that 
  simulating such devices numerically
  (which is necessary for computer-aided design, needed in turn  to
  predict performance and stability prior to manufacturing)
    is performed nowadays in the frequency domain. Indeed, due to the tremendous number of  electronic components
    and high frequency of input signals, such simulations
    can no longer be performed accurately in the time domain, as it
    would require too small a time step. } 
\end{enumerate}


\smallskip

Our result supersedes the ``Henry-Hale theorem'', that settles exponential
  stability issues for difference-delay systems 
  in the \emph{time-invariant case};
  \textit{i.e.},  when the matrices $D_j$ in \eqref{systae} are constant.
This result, stated further below as Theorem \ref{theorem_Hale} in Section \ref{section_results},
dates back to the 1970s; it
was first established in \cite{CRUZ1970} for finitely
many commensurate delays and later carried over to
countably many, not necessarily commensurate delays in
\cite{Henry1974}; generalizations to distributed delay systems may be
found in \cite[Chapter 9]{Hale}.
The stability of functional differential equations like \eqref{eq:2}
  has been studied since the
fifties  in the
time-invariant case, using either Laplace transform and semigroups to reduce
the problem  to localizing the zeros of   almost periodic holomorphic functions
\cite{Bellman1959, Henry1974, CRUZ1970,Hale},
or else devising  Lyapunov–Krasovskii functionals
to prove stability in special cases~\cite{Hale,Mondie,Rasvan}. 
The stability of \emph{time-varying} linear difference-delay
systems was not investigated
nearly as extensively, and the literature we know of can be
broken up as follows\seb{, in three points.}
\emph{First}, a (fairly restrictive) sufficient condition for exponential stability,
based on the Perron-Frobenius theorem, can be found in \cite{Ngoc2015}.
\emph{Second}, a formula representing the  solutions of  \eqref{systae} in the general
time-varying case is given  in
\cite{Chitour2016}.  The latter paper contains
 interesting results on the insensitivity of
$L^p$-exponential stability to $p$, and on how stability is preserved
under perturbation of the delays, but unfortunately the characterization of exponential stability
proposed there is  untractable for it involves the spectrum of sums of
products of matrices with indefinitely growing number of terms and factors
when $t$ goes large, indexed according to combinatorial rules involving
the lattice generated by the real numbers $\tau_j$.
\emph{Third}, if a linear time-varying system of the form \eqref{systae}
arises from a network of lossless telegrapher's equations (for
instance as the high frequency limit of an electrical network),  then a
sufficient condition for exponential stability can be based on
dissipativity at the nodes of the network \cite{BFLP}. Related, somewhat
specific criteria 
for hyperbolic 1-D systems may be found in \cite{CoNg}.


\bigskip

In this paper we dwell on control-theoretic ideas, basing our approach on the introduction
of an associated control system whose exponential stability is equivalent to the one of \eqref{systae}
(adding a virtual control in the simplest possible way, see \eqref{eq:4}).
%
This fresh point of view is suggestive of
  new  tools to investigate stability of periodic systems, as we now explain.
In the time-invariant case, the
transfer function of the associated control system is a matrix-valued function $H(p)$ of a
complex variable $p$, 
and the Henry-Hale theorem amounts to saying that system \eqref{systae} is
exponentially stable if and only if $H$ is holomorphic in some half plane
$\{p\in\corpsC\,,\ \Re(p)>\alpha\}$, with $\alpha<0$;
this is equivalent to the seemingly stronger requirement that $H$ be holomorphic \emph{and bounded} in such a half-plane,
see discussion after Theorem \ref{theorem_hale_generalization}. 
To adress the periodic case,
we define a so-called \emph{Harmonic Transfer Function} (HTF)
for periodic difference-delay control systems
(whose state space is infinite-dimensional),
that generalizes  the one  introduced in \cite{wereley1990analysis}
for periodic differential equations on $\corpsR^n$ (whose state space is finite-dimensional).
Our HTF is a holomorphic map of a complex variable $p$, valued in the
space of linear operators
on $L^2([0,T],\corpsValeurs^d)$ with $T$ the period
  of the system,
that reduces to multiplication by the ordinary transfer function at $p$
when the system is time-invariant.
Then, our main result can be  interpreted in terms of the HTF
the same way as the Henry-Hale theorem does in terms of the
classical transfer function.
Namely, a periodic difference-delay system is 
exponentially stable if and only if the HTF of the associated  control
  system is holomorphic
and bounded (as an operator-valued map) in some half-plane
$\{p\in\corpsC\,,\
\Re(p)>\alpha\}$ with $\alpha<0$. This
  is the content of   Theorem~\ref{theorem_hale_generalization2}, which
  is the main result of the paper; an equivalent  formulation
  not mentioning the HTF (but featuring the latter in disguise)
  is given in Theorem~\ref{theorem_hale_generalization}.


Hence, just like transfer functions encode in their singularities
   the stability properties of time-invariant
  linear systems, harmonic transfer functions
   as defined in this paper   reflect
the stability of infinite-dimensional linear periodic systems.
In this connection,
  we note that results from this paper
 were implicitly anticipated by the Engineering community when
  basing stable design of active circuits on the location of the singularities of
  certain analytic functions that are none but
  the first few Fourier coefficients of the HTF, see  \cite{Suarez} and the bibliography therein.
  In turn, the present study
questions Engineering practice; {\it e.g.}, asking when the singularities of the HTF coincide with those of finitely many such Fourier coefficients.

Besides classical tools from Fourier or complex analysis
  and periodic evolution families,
  the proofs appeal to further material like variation-of-constant formulas in the $BV$-setting  and controlled inversion in Banach algebras\footnote{Controlled inversion (bounding the norm of the inverse in a subalgebra
in terms of the norm in this subalgebra and the norm of the inverse in the algebra) was pioneered in Baskakov \cite{baskakov1997asymptotic} and Nikolski \cite{Nikolski1999}.}, along with
  the fact that $L^p$-exponential stability is independent of $p$
  for systems like \eqref{systae}.
  The most convoluted, and  perhaps deeper part of this work
  is the connection, given by Lemma \ref{prop_lien_mon_ITFbis},
  between the HTF
  and the ordinary transfer function of the lifting of the solution operator
  for \eqref{systae}  (which is a discrete time-invariant infinite-dimensional
  linear system).
  Because the dynamics of (a suitable realization of)
  this lifting is the monodromy operator (see Theorem  \ref{prop_rel_mon_HTF2}), the above-mentioned
  connection is the
  ultimate reason why the singularities of the HTF reflect the spectrum of the monodromy operator and therefore also the exponential
  stability properties of the system.

\bigskip

The paper is organized as follows.
  In Section \ref{section_results}, we make pieces of notation and define
  exponential stability, before recalling the Henry-Hale theorem; we then
  state  its generalization to the periodic case which constitutes our main result
  (Theorem~\ref{theorem_hale_generalization}), before defining the HTF and
  reformulating this main result
  in terms of
  control systems  (Theorem~\ref{theorem_hale_generalization2}).
  Section \ref{sec:notation} contains basic facts on functions of bounded variation, while Section \ref{sec:prelim} introduces fundamental solutions  and
  {\it a priori} estimates thereof, as well as variation-of-constant formulas, for systems of the form \eqref{systae}. The proof of
  Theorems~\ref{theorem_hale_generalization}  and  \ref{theorem_hale_generalization2} is given in
  Section \ref{proofGHeHa}, then Section~\ref{sec:conj} concludes with
a discussion of  a conjecture
regarding neutral periodic delay equations.

\section{Statement of the Main Result}
\label{section_results}

\subsection{Notations}
\label{sec-2-notations}

\newcounter{subsubsecperso}
\setcounter{subsubsecperso}{0}
\renewcommand\thesubsubsecperso{\thesubsection.\arabic{subsubsecperso}}

\noindent
\refstepcounter{subsubsecperso}
\textit{\thesubsubsecperso}.
The real and complex fields are denoted by $\corpsR$ and $\corpsC$. We write $\|\cdot\|$ for Euclidean norm on $\corpsC^d$
and  $\vertiii{\cdot}$ for the norm of an operator or a matrix $\corpsC^d\to\corpsC^\ell$:
$\vertiii{M}=\sup_{\|x\|=1}\|Mx\|$. We put $I_d$ for the $d\times d$ identity matrix or 
operator on a vector space of dimension $d$, and  $I_{\infty}$  for
the ``doubly infinite'' identity matrix
$[\delta_{i,j}]_{(i,j)\in\ZZ^2}$ or the identity
operator on the Hilbert space $\ell^2(\mathbb{Z},\corpsValeurs^d)$.

\medskip

\noindent
\refstepcounter{subsubsecperso}
\textit{\thesubsubsecperso}.
For $E\subset \corpsR$
a Lebesgue-measurable set, {
  we write $L^q(E)$
for the space 
of (equivalence classes of a.e.\ coinciding) $\corpsC$-valued
measurable functions on $E$
with norm
$\|g\|_{L^q(E)}=(\int_E|g(y)|^qdy)^{1/q}$ ($\textrm{ess. sup}_E\,|g|$
if $q=\infty$). 
The space $L^q_{\mathrm{loc}}(E)$
consists of functions whose restriction
to any compact $K\subset E$ lies in $L^q(K)$.
We write $C^0(E)$ for the space 
of $\corpsC$-valued continuous  functions on $E$,
and if $E$ is compact we endow it
with the \emph{sup} norm denoted by $\|\cdot\|_{C^0}$.
For $\alpha\in(0,1)$, we designate with
$C^{\alpha}(E)$ the subspace of  $C^0(E)$ consisting
of H\"older continuous functions
with  exponent $\alpha$; {\it i.e.}, $f\in C^{\alpha}(E)$ if and 
only if
$|f(x)-f(y)|\leq C|x-y|^\alpha$ for some constant $C$ and all $x,y\in E$,
the smallest $C$ being the H\"older constant of $f$.
When $E$ is open, $C^1(E)$ 
indicates the space of complex functions whose
first derivative lies in $C^0(E)$, and
$C^{1,\alpha}(E)$ stands 
for functions whose first derivative belongs to $C^{\alpha}(E)$.
When dealing with
vector-valued functions, we indicate the target space as in
$L^q(E,\corpsValeurs^{d\times d})$ or $C^{1,\alpha}(E,\corpsValeurs^d)$,
while replacing in the definition 
the modulus by the corresponding norm.

\medskip

\noindent
\refstepcounter{subsubsecperso}\label{rappel-l2}
\textit{\thesubsubsecperso}.
We shall work in the Hilbert space:
\begin{eqnarray}
\ell^2(\mathbb{Z},\corpsValeurs^d) := \{ z=(z_j)_{j \in \mathbb{Z}}:\, z_j \in \corpsValeurs^d, \sum\limits_{j=-\infty}^{+\infty}\|z_j\|^2 <+ \infty \}\,,
\end{eqnarray}
equipped with the standard norm
\begin{eqnarray}
\|z\|_{\ell^2}:=  \biggl(\,\sum_{j=-\infty}^{+\infty}\|z_j\|^2\biggr)^{1/2}.
\end{eqnarray}
The norm of bounded operators
$\ell^2(\mathbb{Z},\corpsValeurs^d)\to\ell^2(\mathbb{Z},\corpsValeurs^d)$
(often identified with doubly-infinite $d\times d$ block matrices representing them in a basis)
is denoted by $\vertiii{\cdot}_2$; {\it i.e.}, $\vertiii{L}_2=\sup_{\|z\|_{\ell^2}=1}\|Lz\|_{\ell^2}$.

\medskip

\noindent
\refstepcounter{subsubsecperso}
\textit{\thesubsubsecperso}.
If a function $f$ is defined on $E$ and $E^\prime\subset E$,
we put $f_{|E^\prime}$ to mean the restriction of $f$ to $E^\prime$.

\medskip

\noindent
\refstepcounter{subsubsecperso}\label{rappel-Hardy}
\textit{\thesubsubsecperso}.
The Hardy space $\mathcal{H}^2$ of the right half-plane is comprised  of 
those holomorphic functions $f$ in $\{z\in\corpsC:\,\Re(z)>0\}$
satisfying
\[
\|f\|^2_{\mathcal{H}^2}:=\sup_{x>0}\int_{-\infty}^{+\infty}|f(x+iy)|^2dy<+\infty.
  \]
  Such functions are $\sqrt{2}$-metrically the Laplace transforms of square integrable functions on $[0,+\infty)$ \cite[Ch. 8, p. 131]{Hoffman}. That is, $\mathcal{H}^2=\{f:f(z)=\int_0^{+\infty}e^{-zt}u(t)dt,\,u\in L^2([0,+\infty),\corpsC),\,\Re(z)>0\}$,
  and it holds that $\|f\|^2_{\mathcal{H}^2}=2\|u\|^2_{L^2([0,+\infty),\corpsC)}$.

\subsection{Solution operators and exponential stability}
\label{section_PbStatement}

Consider a  periodic difference-delay system:
\begin{eqnarray}
\label{system_lin_formel}
  y(t)=\sum_{j=1}^ND_j(t)y(t-\tau_j), \qquad t > s,\qquad
y(s+\theta)=\phi(\theta)\  \mathrm{for}\  -\tau_N\leq \theta \leq0,
\end{eqnarray}
where $s\in \corpsR$ is the initial time, $d$ and $N$ are positive 
integers, $\tau_1 < \cdots < \tau_N$  are strictly positive real numbers
(the delays) and
the $D_j:\corpsR\to\corpsValeurs^{d\times d}$
are continuous $T$-periodic matrix-valued functions:
\begin{equation}
  \label{eq:1}
  D_j(t+T)=D_j(t)\,,\ \ 1\leq j\leq N\,,
\end{equation}
while $\phi:[-\tau_N,0]\to\corpsValeurs^d$ is the initial condition
and  solutions to \eqref{system_lin_formel} are
$\corpsValeurs^d$-valued functions $y(t)$  of the time $t\in[s-\tau_N,+\infty)$.
When the $D_j$ are real-valued, all results below  specialize to
real solutions obtained by
restricting to real initial conditions.
We shall assume that $T$ is  strictly larger than the delays:
 $\tau_N<T$,
which is no  loss of generality for we may replace $T$ by $kT$ with $k\in\N$.

One may seek solutions of \eqref{system_lin_formel} in various functional spaces.
As the $D_j(.)$ are  continuous, we may for instance look for continuous solutions in which case 
a compatibility condition is required on the initial condition:
it must
belong to the space
\begin{eqnarray}
C_s :=\{ \phi \in  C^0([-\tau_N,0],\corpsValeurs^d):\,\phi(0)=\sum\limits_{j=1}^ND_j(s)\phi(-\tau_j)\}.
\end{eqnarray}
Given $\phi \in C_s$, an easy recursion shows that system \eqref{system_lin_formel} has 
a unique continuous solution $y\in C^0([s-\tau_N,+\infty),\corpsValeurs^d)$ with initial condition $\phi$. Thus, we can
define for $t\geq s$ the \emph{solution operator}
$U(t,s) : \,C_s \to C_t$,
mapping $\phi\in C_s$ to $U(t,s)\phi$
defined by
\begin{equation}
   \label{defSO}
           \bigl(U(t,s)\phi\bigr)(\theta)= y(t+\theta)\,,
           \quad\theta\in[-\tau_N,0]\,.
\end{equation}
Note that $C_{s+T}=C_s$ and $U(t,s)=U(t+T,s+T)$, by 
the $T$-periodicity of  the $D_j$.

One may also seek solutions  of \eqref{system_lin_formel}
in $L^q_{\mathrm{loc}}([s-\tau_N,+\infty),\corpsValeurs^d)$ for 
$1\leq q\leq\infty$, and then we require
\eqref{system_lin_formel} to hold for almost every $t\in[s,+\infty)$ only;
no compatibility
condition on $\phi$ is needed anymore, and a recursive argument
shows that for each $\phi \in
L^q([-\tau_N,0],\corpsValeurs^d)$ the system \eqref{system_lin_formel} admits  a
unique solution $y\in L^q_{\mathrm{loc}}([s,+\infty),\corpsValeurs^d)$ with initial condition $\phi$.
Consequently,
one can define for $t\geq s$ the solution operator 
$U_q(t,s) : \,L^q([-\tau_N,0],\corpsValeurs^d) \to L^q([-\tau_N,0],\corpsValeurs^d)$
that maps $\phi$ to $U_q(t,s)\phi$
given by a relation  analogous to \eqref{defSO}:
\begin{equation}
  \label{defSOq}
  \bigl(U_q(t,s)\phi\bigr)(\theta)= y(t+\theta)\,,\quad \text{a.e. }\ \theta\in[-\tau_N,0]\,.
\end{equation}
These different types of solutions {\it a priori} yield
  distinct notions of exponential stability defined as follows.
\begin{definition}\label{def:stab}
    System $(\ref{system_lin_formel})$ is called $C^0$-exponentially stable if there exist $\gamma,K>0$ such that 
  \begin{equation}  \label{eq:stabC0}
    \|U(t,s)\phi\|_{C^0} \le Ke^{-\gamma (t-s)} \|\phi\|_{C^0},
    \mbox{ for all $s\in \corpsR$, all $t \ge s$ and all $\phi \in C_s$}.
  \end{equation}
System $(\ref{system_lin_formel})$ is called $L^q$-exponentially stable, $q\in[1,\infty]$, if there exist $\gamma,K>0$ such that 
  \begin{equation}    \label{eq:stabLq}
 \|U_q(t,s)\phi\|_{L^q} \le Ke^{-\gamma (t-s)} \|\phi\|_{L^q}, \mbox{ for all $s\in\corpsR$, all $t \ge s$ and all $\phi \in L^q([-\tau_N,0],\corpsValeurs^d)$}.
\end{equation}
\end{definition}

It is remarkable that these notions are in fact equivalent,
as shown
by the following
result contained in \cite[Theorem 3.4]{BFLP}
which is also a   consequence of  \cite[Corollary 3.29]{Chitour2016}.
\begin{proposition}
\label{equi_L2_C0} 
For each $q\in[1,\infty]$, System $\eqref{system_lin_formel}$ is  $L^q$-exponentially stable if and only it is $C^0$-exponentially stable. 
\end{proposition}

Proposition \ref{equi_L2_C0} plays an important role
in the proof of Theorem~\ref{theorem_hale_generalization},
because the sufficiency part 
establishes $C^0$-exponential stability 
while the necessity part assumes $L^2$-exponential stability.
In view of Proposition \ref{equi_L2_C0},
hereafter we simplify terminology by making  the following definition.

\begin{definition}[Exponential stability]
  \label{defexpss}
System $(\ref{system_lin_formel})$ is called 
\emph{exponentially stable} if and only if if there exist $\gamma,K>0$
such that one of the equivalent properties
\eqref {eq:stabC0} or \eqref{eq:stabLq} holds.
\end{definition}

\subsection{The Henry-Hale Theorem}
In the time-invariant case,
the following characterization of exponential stability is known.
%
\begin{thm}[Henry-Hale {\protect \cite[Section 3]{Henry1974}, \cite[Theorem 3.5]{Hale}}]
  \label{theorem_Hale}
  Assume that the maps $t\mapsto D_j(t)$ are constant, $1\leq j\leq N$.
  Then,
a necessary and sufficient condition for System \eqref{system_lin_formel}
to be $C^0$-exponentially stable is the existence 
of a real number $\beta<0$ such that :
\begin{eqnarray}
\label{cdt0}
  I_d-\sum\limits_{j=1}^N e^{-p \,\tau_j} \,D_j \ \text{is invertible in
  $\corpsValeurs^{d\times d}$ for every $p$ in $\{ z \in \corpsC:\,\Re(z)\ge\beta\}$
  }.
\end{eqnarray}
\end{thm}



Our goal is to ``generalize'' Theorem~\ref{theorem_Hale}
  to the case where the maps $t\mapsto
  D_j(t)$ are \emph{not} constant, but periodic.
Let us first discard the naive attempt requiring that \eqref{cdt0} holds for all $t$, for
this is not enough to ensure
  stability as the following example shows.

Let $N=1$, $d=2$, $T=2$, $\tau=1$, and set $D_1(t)=\left(
\begin{smallmatrix}
  1/2& a(t) \\ b(t) &1/2
\end{smallmatrix}
\right)$
in \eqref{system_lin_formel}, with $a$ a smooth function such that
$a(t)\equiv0$ if $t\in[1,2]$, $a(t)\equiv1$ if $t\in[\frac13,\frac23]$, and
$a(t+2)=a(t)$ for all $t$, while $b(t)=a(t+1)$.
On the one hand, condition \eqref{cdt0} is satisfied because, since
$a(t)b(t)=0$ for all $t$, $D_1(t)$ always has $\frac12$ as a double eigenvalue.
On the other hand, for $t$ in $[\frac13,\frac23]$, one has
$y(t+2k)=\bigl(D_1(t+2)D_1(t+1)\bigr)^ky(t)$ for all $k\in\N$, with
$D_1(t+2)D_1(t+1)=\left(
\begin{smallmatrix}
  5/4& 1/2\\ 1/2 &1/4
\end{smallmatrix}
\right)$; since this matrix
has one eigenvalue larger than 1 (namely $\frac34+\frac{\sqrt{2}}2$), 
exponential stability cannot hold.

In the next section, we shall give a proper analog
  of Theorem~\ref{theorem_Hale} in the periodic case.

\subsection{Main result: generalizing  the Henry-Hale theorem
  to the time-varying periodic case.}
From now on, we consider the case where the maps $t\mapsto D_j(t)$ are
$T$-periodic.
%
Let us define
\begin{equation}
  \label{eq:omega}
  \omega:=2\pi/T\,,
\end{equation}
and put $\check{D}_j(k)$,
for $j\in\{1,\cdots,n\}$ and  $k \in \mathbb{Z}$, to indicate the $k$\textsuperscript{th} Fourier coefficient of $D_j$:
\begin{eqnarray}
\label{fourier_coeef}
\check{D}_j(k):= \frac{1}{T} \int_0^T D_j(t) e^{-i  k \omega t} dt\,.
\end{eqnarray}
We denote by $L_{D_j}$ the (doubly infinite)
block Laurent matrix associated with $D_j(t)$;
its block entries are:
\begin{eqnarray}
  \label{defLDj}
\left(L_{D_j}\right)_{k,\ell}:=\check{D}_j(\ell-k),\qquad{\ell,k \in \mathbb{Z}}.
\end{eqnarray}
We also put $\Delta_{\tau_j,\omega}$, $j\in\{1,\ldots,N\}$, for
the (doubly infinite) block diagonal matrix given by
\begin{eqnarray}
  \label{defDtilde}
\Delta_{\tau_j,\omega}:=\mathrm{\ diag}\bigl(\cdots,\,e^{-2i \omega \tau_j}I_d\,,\,e^{-i \omega \tau_j}I_d\,,\,I_d\,,\,e^{+i \omega \tau_j}I_d\,,\,e^{+2i \omega \tau_j}I_d\,,\cdots\bigr).
\end{eqnarray}
Both $L_{D_j}$ and
$\Delta_{\tau_j,\omega}$ define, by matrix multiplication, bounded
operators
$\ell^2(\mathbb{Z},\corpsValeurs^d)\to\ell^2(\mathbb{Z},\corpsValeurs^d)$.
In fact,
when  $\ell^2(\mathbb{Z},\corpsValeurs^d)$ gets identified with
$L^2([0,T),\mathbb{C}^d)$ (the space of square integrable $\corpsValeurs^d$-valued
functions on the circle of circumference $T$) {\it via} the Fourier coefficients
(arranged columnwise so that indices increase from top to bottom),
then $\Delta_{\tau_j,\omega}$ corresponds to the isometry
$f(\xi)\mapsto f(e^{i\tau_j}\xi)$ while
$L_{D_j}$ becomes pointwise multiplication by $D_j(-t)$, whence
$\vertiii{L_{D_j}}_2=\|D_j\|_{C^0}$. Hereafter, we often identify operators
$\ell^2(\mathbb{Z},\corpsValeurs^d)\to\ell^2(\mathbb{Z},\corpsValeurs^d)$
and operators $L^2([0,T),\mathbb{C}^d)\to L^2([0,T),\mathbb{C}^d)$,
{\it via}  their matrix in the Fourier basis. 



We define a function $R$ of the complex variable $p$,
valued in the space of bounded operators $\ell^2(\mathbb{Z},\corpsValeurs^d)\to \ell^2(\mathbb{Z},\corpsValeurs^d)$  by the formula
\begin{eqnarray}
\label{defHc}
  R(p):=I_{\infty}- \sum\limits_{j=1}^N e^{- p \,\tau_j} L_{D_j} \Delta_{\tau_j,\omega}\,,
\end{eqnarray}
and the doubly infinite matrix  representing $R(p)$ in the Fourier basis
has the block description:
\begin{equation}
  \label{defHc-bis}
  R(p)=
  I_{\infty}-\left(\sum_{j=1}^{N}
    e^{-(p\,-\,i\,\ell\,\omega)\tau_j}
    \check{D}_j(\ell-k)
  \right)_{(k,\ell)\in\mathbb{Z}^2}.
\end{equation}
 We remark \seb{that, by a Neumann series argument, the  infinite matrix $R(p)$ is invertible for all complex number $p$ with large} enough real part,  and we denote by $R(p)^{-1}$ the inverse matrix.
Our main result is now the following.
\begin{thm}[Necessary and sufficient condition for exponential stability]
\label{theorem_hale_generalization}
Assume that the $D_j:\mathbb{R}\to\corpsValeurs^{d\times d}$ are  periodic
with H\"older continuous derivative for $1\leq j\leq N$.
Then, a necessary and sufficient condition for System \eqref{system_lin_formel}
to be 
exponentially stable is the existence 
of a real 
number $\beta<0$ such that :
\begin{enumerate}[label=\textit{(\roman*)}] 

\item \label{assumption1}
$R(p)$ 
is invertible $\ell^2(\mathbb{Z},\corpsValeurs^d) \to \ell^2(\mathbb{Z},\corpsValeurs^d)$ for each $p$ in $\{ z \in \corpsC:\,\Re(z)\ge\beta\}$,

\item \label{assumption2}
there is a positive number $M$ such that 
  $\displaystyle\vertiii{R(p)^{-1}}_2 \leq M$ 
  for all $p$ in $\{ z \in \corpsC:\,\Re(z)\ge\beta\}$.

\end{enumerate}


\end{thm}

When the $D_j$ are constant matrices,
the operators $L_{D_j}\Delta_{\tau_j,\omega}$ are block diagonal
with $k$\textsuperscript{th} diagonal block $e^{ik\omega\tau_j}D_j$ and condition $(i)$ in
Theorem~\ref{theorem_hale_generalization} is equivalent to
saying that  $I_d-\sum_{j=1}^N e^{-p\tau_j}D_j$ is invertible for $\Re(p)\geq\beta$.
In this case, since a holomorphic function on a vertical strip
which is uniformly almost-periodic in the pure
imaginary direction and has no zeros must be bounded below in modulus
by a strictly positive constant \cite[Ch. III, Sec. 2, Cor. I]{besicovitch1954almost},  it holds that $|\mathrm{det}(I_d-\sum_{j=1}^N e^{-p\tau_j}D_j)|\geq c>0$
for $\Re(p)\geq\beta$, whence
condition $(ii)$
is redundant. Thus, 
Theorem~\ref{theorem_hale_generalization}  yields back the Henry-Hale theorem
for systems with constant coefficients.
\seb{More generally, Condition \ref{assumption1} implies
Condition \ref{assumption2} 
when all delays are commensurate.} Indeed, $p \mapsto R(p)$
is then periodic of period $i\tau$ for some $\tau>0$ and hence, since inversion
is continuous on invertible elements of a Banach algebra while
\ref{assumption2} needs only be checked for $p$ in a
compact set  by periodicity and a Neumann series argument, 
we deduce that \ref{assumption1} implies \ref{assumption2} in this case.
Thus, we obtain the following:

\begin{corollary}
  If the delays are commensurate, then
  condition \ref{assumption2} can be omitted
  in Theorem~\ref{theorem_hale_generalization}.
\end{corollary}

For  periodic \seb{non-commensurate} difference-delay systems, the authors doubt that Condition \ref{assumption1} implies
Condition \ref{assumption2} in general, though they know of no counterexample.
As mentioned above, in the time-invariant case the redundancy of $(ii)$
  comes from \seb{properties of  almost-periodic (complex valued) holomorphic}
  functions,
  \seb{for which} the values at infinity are
  linked to the values at finite distance thanks to \seb{almost periodicity
  combined with}  theorems of Montel and Rouché.
  Unfortunately, \seb{no straightforward extension to almost-periodic \emph{operator valued} complex analytic
    functions 
  is available in general, because both the Montel and Rouché theorems fail, at least     in a non-Fredholm context.}

\begin{rmk}
  The proof of Theorem~\ref{theorem_hale_generalization} will show that
  the sufficiency part
  remains true when the $D_j$ are merely continuous and the  assumption \ref{assumption2} is replaced by:
  \begin{eqnarray}
\label{eq:th2_norm_Wiener}
    \vertiii{
 R(p)^{-1}}_{\mathfrak{W}} \mbox{ is uniformly bounded for all $p$ in } \{ z \in \corpsC:\,\Re(z)\ge\beta\},
 \end{eqnarray}
 where the Wiener norm $\vertiii{\cdot}_{\mathfrak{W}} $ of a doubly infinite block matrix is  defined
 in \eqref{def_norm_Wiener}.
\end{rmk}

\subsection{A system-theoretic point of view; transfer functions and
  harmonic transfer functions.}
\label{sec:control}

  Below we recast the previous considerations  in} system-theoretic language,
introducing harmonic transfer functions and reformulating
Theorem~\ref{theorem_hale_generalization} as
Theorem~\ref{theorem_hale_generalization2}, which is more faithful to the
version of our main
result described in the introduction.

To  \eqref{systae} one can associate the control system:
\begin{equation}
     \label{eq:4}
     y(t)=\sum\limits_{j=1}^ND_j(t)y(t-\tau_j)+u(t),\qquad  t > s,
     \qquad y(t)= u(t)=0\ \mathrm{for}\ t< s,
\end{equation}
with control $u \in C^0([s,+\infty),\corpsValeurs^d)$ (or $u\in L^p_{loc}([s,+\infty),\corpsValeurs^d)$) and output $y$. System \eqref{eq:4} is
 exponentially stable if, when driven
by an input $u(t)$ vanishing for $t>t_0$ and generating an output $y$,
the $L^\infty$ (or $L^q$)-norm of the restriction $y_{|[t-\tau_N,t]}$
decays exponentially fast  to zero as $t\to+\infty$.

Solving for $u$ in \eqref{eq:4},
   one sees that for any  $t_0>\tau_N$ and  $\phi\in  C_{t_0}$
    (resp.\ $L^p([-\tau_N,0],\corpsValeurs^d)$),
   there is $u\in C^0([s,t_0],\corpsValeurs^d)$
    (resp.\ $L^p([s,t_0],\corpsValeurs^d)$)
   such that the corresponding output $y$ satisfies 
    $\phi(\theta) =y(t_0+\theta)$ for $\theta\in[-\tau_N,0]$.
   Consequently,  outputs of \eqref{eq:4} associated to controls that vanish for $t>t_0$ coincide with solutions to \eqref{systae} where we put $s=t_0$.
   This explains why  exponential stability of the control system
   \eqref{eq:4} is equivalent to  exponential stability of  \eqref{systae}.

A linear time-invariant control system
can be represented, under mild assumptions, as a convolution operator;
{\it i.e.}, the output is obtained by convolving the input with some kernel
\cite{RudinFA}. Taking Fourier-Laplace transforms converts the latter into a multiplication operator in the Fourier-Laplace domain, and then
  the multiplier is called the \emph{transfer function} of the system.
Hereafter we specialize this to time-invariant delay difference control systems, before explaining
  the notion corresponding to the transfer-function in the 
  periodic case.

\emph{In the
  time-invariant case}, the transfer function of \eqref{eq:4} is  the matrix-valued function of one complex
variable, say $p$, given by 
$H(p)=\left(I_d-\sum_{j=1}^N e^{-p \,\tau_j} \, D_j\right)^{-1}$
(see \textit{e.g.} \cite{Faur-Dep76}).
The Henry-Hale theorem
says that \eqref{system_lin_formel} is exponentially stable if and only if
$H$ is holomorphic in a half-plane $\{p\in\corpsC\,,\
\Re(p)>\beta\}$ for some $\beta<0$. In this case, as explained after
the statement of Theorem \ref{theorem_hale_generalization}, $H(p)$ is in fact bounded
for $\Re(p)\geq\beta'>\beta$.

\emph{In the periodic case}, the concept of transfer function
generalizes into the one of  \emph{harmonic transfer function} (in short: HTF).
It is again a function of one complex variable $p$, but instead of ranging
in $\corpsValeurs^{d\times d}$ it takes values in the space of bounded operators
$\ell^2(\ZZ,\corpsValeurs^d)\to \ell^2(\ZZ,\corpsValeurs^d)$.
More precisely,  one  can check from \eqref{eq:4}  that
$\|y\|_{L^2([s,\tau])}\leq Ce^{\gamma (\tau-s)}\|u\|_{L^2([s,+\infty])}$ for
appropriate constants $C$, $\gamma$ and all $\tau>s$;
see \eqref{majexpL2py}.
Thus, when  $u \in L^2([s,+\infty),\corpsValeurs^d)$ one can define,
for $\Re(p)>\gamma$,  the Laplace transforms of $y(t)$ and $u(t)$ by the formulas
\begin{eqnarray}
\hat{Y}(p):=\int_{s}^{+ \infty} e^{-pt} y(t)dt \quad\mbox{ and }\quad
\hat{U}(p):=\int_{s}^{+ \infty} e^{-pt} u(t)dt
\end{eqnarray}
(recall from \eqref{eq:4} that $u,y$ vanish on $(-\infty,s]$).
Sampling $\hat{Y}$, $\hat{U}$ 
at equally spaced points on vertical lines with mesh $i\omega$
(recall \eqref{eq:omega}), we  construct infinite colums  of
$d\times 1$ vectors:
\begin{eqnarray}
\label{ltvs}
 \hat{\mathbb{Y}}(p):=\left( \begin{array}{c}
\vdots \\
\hat{Y}(p+ i\omega)\\
\hat{Y}(p)\\
\hat{Y}(p-i\omega) \\
\vdots 
                          \end{array} \right)\qquad\text{and}\qquad
   \hat{\mathbb{U}}(p):=\left( \begin{array}{c}
\vdots \\
\hat{U}(p+i\omega)\\
\hat{U}(p)\\
\hat{U}(p-i\omega) \\
\vdots 
\end{array} \right),\qquad \Re(p) >\gamma.
\end{eqnarray}
We choose to order the entries of doubly infinite vectors  from top to bottom;
\textit{i.e.}, we write
\begin{equation}
  \label{eq:6}
  \hat{\mathbb{Y}}(p)=\left( \hat{Y}(p-i k \omega)\right)^{\mathbf{t}}_{k\in\mathbb{Z}}
  \quad\mathrm{and}\quad 
  \hat{\mathbb{U}}(p)=\left( \hat{U}(p-i k \omega)\right)^{\mathbf{t}}_{k\in\mathbb{Z}},
\end{equation}
where
the superscript $\bf t$ means transpose.

Since $\hat{U}(p+ik\omega)/T$
(resp.\ $\hat{Y}(p+ik \omega)/T$) is the $k$\textsuperscript{th}
Fourier coefficient of the $T$-periodic function
$x\mapsto\sum_{j=-\infty}^{+\infty}e^{-p(x+jT)}u(x+jT)$ (resp.\ $x\mapsto\sum_{j=-\infty}^{+\infty}e^{-p(x+jT)}y(x+jT)$) that lies in
$L^2([0,T])$ for $\Re(p) >\gamma$,
both $\hat{\mathbb{Y}}$ and  $\hat{\mathbb{U}}$
lie in $\ell^2(\ZZ,\corpsValeurs^d)$.
If moreover $\Re(p)>\log(\sum_{j=1}^N\|D_j\|_{C^0})/\tau_1$,
there is an operator-valued holomorphic function
$p\mapsto\mathbfsf{H}(p)$ such that $\mathbfsf{H}(p):\ell^2(\mathbb{Z},\corpsValeurs^d)\to \ell^2(\mathbb{Z},\corpsValeurs^d)$ maps $\hat{\mathbb{Y}}$ to  $\hat{\mathbb{U}}$.
Indeed, $\mathbfsf{H}(p)$ is none but
the inverse of $R(p)$ defined in  \eqref{defHc}, computed {\it via}
  a Neumann series; {\it cf.} Section~\ref{sec:necessariness},
Equation \eqref{syst_ref0}.
Altogether, increasing $\gamma$ if necessary, we get that
\begin{eqnarray}
\label{link_fond_HTF}
 \hat{\mathbb{Y}}(p) =\mathbfsf{H}(p)\, \hat{\mathbb{U}}(p),\quad
  \Re(p) >\gamma\,.
\end{eqnarray}

\begin{definition}[Harmonic Transfer Function]
\label{def:HTF}
The operator-valued holomorphic function
$\mathbfsf{H}(\cdot)$ is called the harmonic transfer function (HTF)
of system \eqref{eq:4}.
\end{definition}
%
%
%
We can now  restate Theorem~\ref{theorem_hale_generalization}
as follows.

\begin{thm}[Reformulation of Theorem~\ref{theorem_hale_generalization}]
  \label{theorem_hale_generalization2}
  If system (\ref{system_lin_formel}) has  periodic coefficients of class $C^{1,\alpha}$, it  is exponentially stable if and only if
  the   control system  \eqref{eq:4} is in turn exponentially stable, and that is  if and only if
the harmonic transfer function of the latter 
is analytic and bounded in $\{p\in\corpsC:\,\,\Re(p)>\beta\}$ for some
  $\beta<0$.
\end{thm}
\begin{pf}
   Equivalence between the exponential stability of \eqref{system_lin_formel}
    and \eqref{eq:4} was explained after
    Equation~\eqref{eq:4}. The remaining assertion is formally equivalent to
    Theorem~\ref{theorem_hale_generalization}, granted the definition of the HTF \seb{and the fact that $p\mapsto R(p)^{-1}$ is analytic on its open domain of definition. Indeed, $A\to A^{-1}$ is analytic on the open set of invertible operators  $\ell^2(\mathbb{Z},\corpsValeurs^d)\to \ell^2(\mathbb{Z},\corpsValeurs^d)$,
      because if $A_0$ is invertible then so is $A_0+\delta A$ for
      $\vertiii{\delta A}_{2} \,\vertiii{A_0^{-1}}_{2}<1$ with inverse
      $(A_0+\delta A)^{-1}=A_0^{-1}\sum_{j=0}^\infty(A_0^{-1}\delta A)^j$; since
    $p\mapsto R(p)$ is clearly analytic as well, so is $p\to
    R(p)^{-1}$ by composition, as claimed.}
    \qed
\end{pf}

A couple of remarks are  in order. First, the harmonic transfer function
  can be defined in the same manner for more general linear periodic systems
than difference-delay ones, but we shall stick to the present setting. Second,
Equation~\eqref{link_fond_HTF} 
connects  inputs and outputs of a linear
periodic control system  in the frequency
domain {\it via} the HTF, and if the system is time-invariant,
then the HTF reduces to the block diagonal matrix
${\bf H}(p)=\mathrm{diag}\{\cdots,H(p+i\omega),H(p),H(p-i\omega),\cdots\}$
where $H$ is the ordinary  transfer function,
so that \eqref{link_fond_HTF}
is equivalent to the well-known  input-output relation 
 $\hat{Y}(p)=H(p)\hat{U}(p)$ for time-invariant systems.
Third, the HTF  further generalizes the ordinary transfer function in a way which is worth explaining:
if a stable, \emph{time-invariant}
control system  is fed with a periodic input signal $e^{i\nu t} v$ with $v\in\corpsValeurs^d$, then the output is asymptotically
$H(i\nu)\,v \,e^{i\nu t}$, where $H$ is the transfer function.
If a stable \emph{periodic}
control system
is fed with a periodic input signal $v\,e^{i\nu t}$, then the output is asymptotically
$\mathbfsf{H}(i\nu)({\bf v})e^{i\nu t}$, where $\mathbfsf{H}$ is the harmonic transfer
function and $\bf{v}$ the constant function with value
$v$ in $L^2([0,T),\mathbb{C}^d)$ (so that $\mathbfsf{H}(i\nu)({\bf v})$ is again an element of $L^2([0,T),\mathbb{C}^d)$).
Thus, while in the time-invariant case a periodic input is asymptotically mapped to a periodic output
with the same period, it gets
asymptotically mapped in the periodic
case to an oscillating signal with the same period but  carried
by the wave $\mathbfsf{H}(i\nu)(\bf{v})$ which
has the period of the system. This attractive interpretation
follows from \eqref{idGH},
\eqref{fourier_serie_ITF}, \eqref{eq_ITF} and \eqref{inp-output_equation}.

\section{Functions with bounded variation and
  Lebesgue-Stieltjes integrals} 
\label{sec:notation}

Hereafter we recall basic facts regarding functions of bounded variation
and Stieltjes integrals, setting up some notation regarding bounds of
integration that will be of use throughout.

For $I$ a bounded real interval and $f: I\to \mathbb{R}$ a function,
the \emph{total variation} of $f$ on $I$ is defined as
\begin{equation}
\label{defvar}
W_I(f):=\sup_{\stackrel{x_0<x_1<\cdots<x_N}{x_i\in I, N\in\mathbb{N}}}\sum_{i=1}^N|f(x_i)-f(x_{i-1})|<\infty.
\end{equation}
The space $BV(I)$ of functions with \emph{bounded variation} on $I$
consists of those $f$ such that  $W_I(f)<\infty$,
endowed with the  norm $\|f\|_{BV(I)}=W_{I}(f)+|f(d)|$ where $d\in I$ is arbitrary but fixed. 
Different $d$ give rise to equivalent norms for which
$BV(I)$ is a Banach space, and 
$\|.\|_{BV(I)}$ is stronger than the uniform norm.
We let $BV_r(I)$ and  $BV_l(I)$) be the closed subspaces
of $BV(I)$ comprised of right and left continuous functions, respectively.
We write $BV_{loc}(\mathbb{R})$ for the space of functions whose restriction to any bounded interval $I\subset \mathbb{R}$ lies in $BV(I)$.
Observe that
\begin{equation}
\label{varprod}
W_I(fg)\leq W_I(f)\sup_{x\in I}|g(x)|+W_I(g)\sup_{x\in I}|f(x)|.
\end{equation}
Each $f\in BV(I)$ has a limit $f(x^-)$ (resp.\ $f(x^+)$) from the left (resp.\ right) at every $x\in I$
where the limit applies \cite[sec. 1.4]{lojasiewicz1988introduction}. Hence,
one can associate to $f$ a
finite signed Borel measure $\nu_f$ on $I$ such that
$\nu_f((a,b))=f(b^-)-f(a^+)$, and if $I$ is bounded on the right (resp.\ left) and contains its endpoint $b$ (resp.\ $a$), then $\nu_f(\{b\})=f(b)-f(b^-)$
(resp.\ $\nu_f(\{a\})=f(a^+)-f(a)$) \cite[ch. 7, pp. 185--189]{lojasiewicz1988introduction}.
Note  that different $f$ may generate the same $\nu_f$: for example 
 if $f$ and $f_1$ coincide except  at isolated interior points
 of $I$,  then $\nu_f=\nu_{f_1}$. For $g:I\to\mathbb{R}$ a measurable 
function,  summable against $\nu_f$, 
the Lebesgue-Stieltjes integral $\int gdf$ is defined as
$\int gd\nu_f$, whence the differential element $df$ identifies with $d\nu_f$
\cite[ch. 7, pp. 190--191]{lojasiewicz1988introduction}. This type of integral
is
useful  to integrate by parts, but  caution
must be used when integrating a function with respect to $df$ over
a subinterval $J\subset I$ because
$\nu_{(f_{|J})}$ needs \emph{not} coincide with the restriction $({\nu_f})_{|J}$ of $\nu_{f}$ to $J$. More precisely, if the lower bound $a$ (resp.\ the upper bound $b$)
of $J$ belongs to $J$ and lies interior to $I$, then the two measures may differ by the weight they put on $\{a\}$ (resp.\ $\{b\}$), and they agree only when
$f$ is left (resp.\ right) continuous at $a$ (resp.\ $b$).
By $\int_J gdf$, we always mean that we integrate $g$ against $\nu_{(f_{|J})}$
and \emph{not} against  $({\nu_f})_{|J}$.
We often trade the notation  $\int_J gdf$ for  one of the form
$\int_{a^\pm}^{b^\pm}gdf$, where  the interval of integration $J$ is encoded
in  the  bounds put on the  integral sign: a lower bound $a^-$ (resp.\ $a^+$) means that
$J$ contains (resp.\ does not contain) its lower bound $a$, while an
upper bound $b^+$ (resp.\ $b^-$) means that
$J$ contains (resp.\ does not contain) its upper bound $b$.
Then, the previous word of caution applies to additive rules: for example,
when splitting $\int_{a^\pm}^{b^\pm}gdf$ into
  $\int_{a^\pm}^{c^\pm}gdf+\int_{c^\pm}^{b^\pm}gdf$ where $c\in(a,b)$, we must use $c^+$ (resp.\ $c^-$) if $f$ is right (resp. left) continuous at $c$.

  To a finite, signed or complex Borel measure $\mu$ on $I$,
  one can associate its \emph{total variation measure} $|\mu|$,
  defined on a Borel set $B\subset I$ by
$|\mu|(B)=\sup_{\mathcal{P}}\sum_{E\in\mathcal{P}}|\mu(E)|$ where $\mathcal{P}$ ranges over all partitions of $B$ into Borel sets, see \cite[sec. 6.1]{Rudin};
its total mass $|\mu|(I)$ is called the total variation of $\mu$,
denoted as $\|\mu\|$.
Thus, the total variation is defined both for functions of bounded variation
and for measures, with different meanings. When $f\in BV(I)$ is
monotonic then $W_{I}(f)=\|\nu_f\|$, but in general it only holds  that
 $\|\nu_f\|\leq2W_{I}(f)$; this follows from the Jordan decomposition of $f$
 as a difference of two increasing functions, each of which has variation at most
 $W_{I}(f)$ on $I$ \cite[Thm.\ 1.4.1]{lojasiewicz1988introduction}. 
 In any case, it holds that
  $|\int gdf|\leq \int|g|d|\nu_f|\leq 2W_I(f)\sup_I|g|$.
  The previous notation and definitions   also apply to 
vector and matrix-valued $BV$-functions, replacing absolute values
in \eqref{defvar} by  Euclidean and operator norms, respectively.

  \section{Variation-of-constants formulas and a priori
    estimates}
\label{sec:prelim}

Below, we introduce fundamental solutions and variation-of-constants
  formulas
  for system \eqref{systae}, before deriving estimates thereof; this material
  will be of much use later in the paper, and
we could not see it proven in the literature.

\subsection{Fundamental solution and variation-of-constants formula for continuous solutions.}

The  \emph{fundamental solution}
$X:\corpsR^2\to\corpsValeurs^{d\times d}$ of \eqref{system_lin_formel}
is defined by the following equation
\cite[Equation~(1.3)]{note-formule} (or
  \cite[Theorem~1.2, Chapter 9]{Hale} adressing more general dynamical
  systems):
 \begin{eqnarray}  
 \label{solution_fondamentale}
X(t,s)=\left\{
 \begin{array}{ll}
 0 \mbox{ for $t<s$}, \\
I_d+\sum\limits_{j=1}^N D_j(t) X(t-\tau_j,s) \mbox{ for $t \ge s$}.
\end{array}
\right.
\end{eqnarray}	
Argueing inductively, it is easy to check that $X$ uniquely exists.
Moreover, it does not grow faster than exponential with respect to $t-s$, as the   following proposition shows.
\begin{proposition}
  \label{prop:fondamental-sousexp}
  There exist  $K>0$ and $\lambda\in\corpsR$ such that
\begin{eqnarray}
\label{bound_sol_fund}
\vertiii{X(t,s)} \leq K e^{\lambda(t-s)}, \mbox{ for all $t \ge s$}.
\end{eqnarray}
\end{proposition}
\begin{pf}
  Since the maps $D_i$ are continuous and periodic, there is a $K'>0$ such that
  $\vertiii{D_i(t)}\leq K'$ for  $i\in\{1,\ldots,N\}$ and $t\in\corpsR$.
  Pick $K\geq 2$ and $\lambda$ large enough that 
  \begin{eqnarray}
    \label{bound_sol_fund_preuve1}
    K'Ne^{-\lambda \tau_1}<1/2. 
  \end{eqnarray}
  Let us prove by induction that, for any $k\in\N$, Equation \eqref{bound_sol_fund} holds
  if $t-s < k \tau_1$; this claim clearly implies what we want. Now, this
  is obvious for $k=0$, because $t-s<0$ whence $X(t,s)=0$ in this case.
  Next, assuming that the claim holds for a certain $k\geq0$ and considering $t,s$
  such that
  $k\tau_1 \le t-s < (k+1) \tau_1$, we get from \eqref{solution_fondamentale} that
  \begin{displaymath}
    \vertiii{X(t,s)} \leq 1+ K'NKe^{-\lambda \tau_1}e^{\lambda(t-s)}\,,
  \end{displaymath}
  which implies \eqref{bound_sol_fund} in view of
  \eqref{bound_sol_fund_preuve1} and the inequality 
  $1\leq \frac12 K e^{\lambda(t-s)}$.
  This concludes the induction.
  \qed
\end{pf}
By inspection, $X$ is as smooth as the maps $D_j(.)$.
In particular, it is continuous
about each $(t,s)$ such that $t-s \notin \mathcal{F}$, where
$\mathcal{F}$ is the positive lattice generated by the $\tau_\ell$ in $\corpsR$:
\begin{eqnarray}
\label{set_discontinuities}
\mathcal{F}:= \bigl\{ \sum\limits_{\ell=1}^N n_\ell \, \tau_\ell \,,\; (n_1,\ldots,n_N)\in\mathbb{N}^N\bigr\}\,.
\end{eqnarray}
Clearly, $X$ has a  bounded jump across each line $t-s=\mathfrak{f}$ for
$\mathfrak{f}\in\mathcal{F}$; in fact, a moment's thinking will convince the reader that it is of the form
    \begin{equation}
    \label{eq:5}
    X(t,s)
    \,=\, - \sum_{\mathfrak{f}\in\,[0,t-s]\,\cap\,\mathcal{F}}\,\mathfrak{C}_\mathfrak{f}(t)\,,
    \quad s\leq t\,,
  \end{equation}
where each $\mathfrak{C}_\mathfrak{f}(.)$ 
is differentiable with  Hölder continuous derivative of the same exponent as the maps
$D_j(\cdot)$.
One can see also that  $\mathfrak{C}_\mathfrak{f}(t)$ is a finite sum of products of matrix-valued functions of the form
  $D_j(t-\mathfrak{f}')$,
where $\mathfrak{f}'$ ranges over the elements of $\mathcal{F}$ whose
defining integers $n_\ell$ in \eqref{set_discontinuities} do not exceed
those defining $\mathfrak{f}$, the
empty product being the identity matrix.
A precise expression for $\mathfrak{C}_\mathfrak{f}$ can be obtained by
reasoning as in \cite[Sec. 3.2]{Chitour2016} or \cite[Sec. 4.5]{BFLP}, but we will not need it.
Note that, for fixed $t$, the function $X(t,\cdot)$ is locally
of bounded variation  on $\corpsR$. More precisely, observe from \eqref{solution_fondamentale} that  
$\alpha\mapsto X(t,\alpha)$ lies in $BV_{loc}(\corpsR)$ for all $t$, and that for
$b>s$ we have:
\begin{eqnarray}  
 \label{solution_fondamentalev}
d_\alpha X(t,\alpha)=\left\{
 \begin{array}{ll}
\sum\limits_{j=1}^N D_j(t) d_\alpha X(t-\tau_j,\alpha)  \mbox{ if $t \ge b$} \\
-I_d\delta_t+\sum\limits_{j=1}^N D_j(t) d_\alpha X(t-\tau_j,\alpha) \mbox{ if $s\leq t <b$}
\end{array}
  \right.
  \quad\mathrm{on}\ [s,b],
\end{eqnarray}	
where $\delta_t$ is the Dirac mass at $t$. 
The fundamental solution is quite important because
it yields explicit integral formulas to express
  solutions of \eqref{system_lin_formel}.
  For instance, Equation \eqref{formula_representation} below
allows one  to parametrize continuous solutions, and
is similar in spirit to variation-of-constants
formulas for autonomous time-invariant linear difference-delay
systems given, say  in \cite{Henry1974} or \cite{Haleone,Hale}. Proving it here would make the paper unbalanced;  however, we  provide an argument
in the separate note~\cite{note-formule}. We do so
because, to the difference of \cite{Henry1974, Haleone,Hale},
we deal with time-varying matrices $D_j(.)$, and also because
many such formulas in the literature seem to have issues.
The integral $\int_{s^-}^{(s+\tau_j)^-}$ in Equation~\eqref{formula_representation}
must be understood as a Lebesgue-Stieltjes integral on the interval $[s,s+\tau_j)$;
{\it cf.} Section~\ref{sec:notation}. This type of integral is especially adequate in the present setting of jump singularities, as it cleanly accounts
for excluding or including endpoints of intervals that
may, or may not be charged by the integrand;
everything  is well defined here because $X(t,.)$ is in $BV_{loc}(\corpsR)$.

\begin{proposition}\rm{{\bf (}\cite{note-formule}{\bf )}}
  \label{prop:fondamental}
  For  $s\in\corpsR$ and  $\phi$ in $C_s$ the 
  solution $y\in C^0([s-\tau_N,+\infty),{\corpsValeurs^{d}})$
  to \eqref{system_lin_formel} is
  \begin{eqnarray}
  \label{formula_representation}
y(t)=-\sum\limits_{j=1}^N \int_{s^{-}}^{(s+\tau_j)^{-}} d_{\alpha}X(t,\alpha)D_j(\alpha)\phi(\alpha-\tau_j-s),\qquad t \ge s\,,
  \end{eqnarray}
  where $X$ was defined in
  \eqref{solution_fondamentale}.
\end{proposition}
%
%

\subsection{Variation-of-constants formula for $L^2_{loc}$ solutions.}

Proving necessity in Theorem \ref{theorem_hale_generalization}
requires that we deal with  $L^2([s-\tau_N,s])$-initial data   as well as
$L^2_{loc}([s,\infty))$-solutions for \eqref{system_lin_formel}, and then the variation-of-constants formula features vector-valued integration.
Specifically, let us consider the control system
with input $u$: 
   \begin{eqnarray}  
  \label{syst_lin_avec_entree0}
y(t)=\sum\limits_{j=1}^ND_j(t)y(t-\tau_j)+u(t), \qquad u(t)\in\corpsValeurs^d,\quad t \ge s.
\end{eqnarray}
When $t \ge s$, we set $y_t(\theta)=y(t+\theta)$ for $-\tau_N\leq \theta\leq0$.
The function $y_t:[-\tau_N,0]\to\corpsValeurs^d$ (more accurately: its equivalence class modulo coincidence almost everywhere) is the state at time $t$
of the input-output system \eqref{syst_lin_avec_entree0}.
In particular $y_s$ is now the initial condition, previously denoted by $\phi$
 when dealing with the homogeneous system  \eqref{system_lin_formel} (for which $u\equiv0$). 

 Recalling $X$ from \eqref{solution_fondamentale}, we
 define for each $(t,\alpha)\in\corpsR^2$ a map
$K(t,\alpha):[-\tau_N,0]\to \corpsValeurs^{d\times d}$ by 
\begin{eqnarray}
\label{noyau_convolutionm}
K(t,\alpha)(\theta)=-X(t+\theta,\alpha),\qquad \theta\in [-\tau_N,0].
\end{eqnarray}
For $t\in\corpsR$, $\theta\in[-\tau_N,0]$ and $I\subset\corpsR$ an interval bounded on the left,
the map $\alpha\mapsto X(t+\theta,\alpha)$ lies in
$BV_l(I)$ and vanishes for $\alpha>t+\theta$. As indicated
in Section \ref{sec:notation}, we  associate to this map
a $\corpsValeurs^{d\times d}$-valued Borel measure $\nu_{(X(t+\theta,\,.\,)_{|I}) }$ on $I$.
It follows from \eqref{solution_fondamentale}, \eqref{eq:5}
 and \eqref{solution_fondamentalev}
that  for fixed $t$ and $\theta$, this measure is of the form
\begin{equation}
  \label{formmeas}
\nu_{(X(t+\theta,\,.\,)_{|I}) }=\sum_{\mathfrak{f}\in\mathcal{F},\, t+\theta-\mathfrak{f}\in \widetilde{I}} \mathfrak{C}_\mathfrak{f}(t+\theta)\delta_{t+\theta-\mathfrak{f}},
\end{equation}
where $\widetilde{I}$ denotes $I$ deprived from its right endpoint
(if contained in $I$, otherwise  $\widetilde{I}=I$)
and $\delta_{t+\theta-\mathfrak{f}}$ is a Dirac mass at
$t+\theta-\mathfrak{f}$.
  For fixed $t$, the boundedness of $I$ from the left implies that
  the number of terms in \eqref{formmeas}
  is majorized independently of $\theta\in[-\tau_N,0]$, and clearly
the map sending a Borel set $E\subset I$ to
the function $\theta\mapsto-\nu_{(X(t+\theta,\,.\,)_{|I}) }(E)$ is a vector measure,
valued in the space $\mathcal{B}([-\tau_N,0],\corpsValeurs^{d\times d})$ of bounded measurable $\corpsValeurs^{d\times d}$-valued
functions on $[-\tau_N,0]$ endowed with the sup norm. In view of
\eqref{noyau_convolutionm} we denote this measure with
$d_\alpha K(t,\alpha)$, and though it depends on 
$I$ the latter will be understood from the context.
Then, for $g\in C^0(I,\corpsValeurs^d)$, one can  define the integral
$\int_IdK(t,\alpha)g(\alpha)$ as a member of $\mathcal{B}([-\tau_N,0],\corpsValeurs^d)$;
see, {\it e.g.}  \cite{Dinculeanu} for a definition of  integrals against a
vector measure. To us, it is just a compact notation for the function
\begin{equation}
  \label{intvv}
  \left(\int_I dK(t,\alpha)g(\alpha)\right)(\theta):=-\int_Id_\alpha X(t+\theta,\alpha)g(\alpha)=-\sum_{\mathfrak{f}\in\mathcal{F},\, t+\theta-\mathfrak{f}\in \widetilde{I}} \mathfrak{C}_\mathfrak{f}(t+\theta) g(t+\theta-\mathfrak{f}). 
  \end{equation}
  The rightmost term in \eqref{intvv} can be rewritten as $\sum_{\mathfrak{f}\in\mathcal{F}}\mathfrak{C}_\mathfrak{f}(t+\theta)\widetilde{g}(t+\theta-\mathfrak{f})$ where $\widetilde{g}$ is equal to $g$ on $\widetilde{I}$ and to zero elsewhere.
  Since the $\mathfrak{C}_\mathfrak{f}$ are bounded (being periodic in $C^{1,\delta}(\mathbb{R},\corpsValeurs^{d\times d})$)
    while  those
  $\mathfrak{f}\in\mathcal{F}$  for which $t+\theta-\mathfrak{f}\in\widetilde{I}$ for some
  $\theta\in[-\tau_N,0]$ are finite in number
  because $I$ is bounded on the left,
  one can check that
  $\|\int_IdK(t,\alpha)g(\alpha)\|_{L^2(I,\corpsValeurs^{d})}\leq C\|g\|_{L^2(I,\corpsValeurs^d)}$ for some constant $C$ depending on $t$,   $I$, the $D_j$ and
  the $\tau_j$. Hence, \eqref{intvv} makes good sense for $g\in L^2(I,\corpsValeurs^d)$, but  of course $(\int_I dK(t,\alpha)g(\alpha))(\theta)$ is only defined for a.e.\ $\theta\in[-\tau_N,0]$ in this case. As
  indicated in Section \ref{sec:notation}, the interval of integration
  will  be encoded in the bounds put on the integral sign.
  Note that it is immaterial here whether $I$ contains its right endpoint or not
  when writing $\int_I dK(t,\alpha)g(\alpha)$,  for
  the rightmost term in \eqref{intvv} depends only on $\widetilde{I}$.
  Moreover, $\int_I dK(t,\alpha)g(\alpha)$ is independent of $g(\alpha)$ for $\alpha>t$ as  $K(t,\alpha)\equiv0$ when $\alpha>t$. Hence,
  if $g\in L^2_{loc}([a,+\infty))$ then
  $\int_{a^-}^{b^\pm} dK(t,\alpha)g(\alpha)$ is independent of $b>t$
  and of the choice of sign in $b^\pm$. In this case, we find it convenient to define
  \begin{equation}
    \label{defborninf}
    \int_{a^-}^{+\infty} dK(t,\alpha)g(\alpha):=
    \int_{a^-}^{b^\pm} dK(t,\alpha)g(\alpha)\qquad
    \text{for any $b>t$.}
    \end{equation}
  
    The variation-of-constants formula for $L^2_{loc}$-solutions to \eqref{syst_lin_avec_entree0}
    now goes as follows, with $U_2(.,.)$  the solution operator defined in~\eqref{defSOq}).
\begin{proposition}
  For $s \in \mathbb{R}$, let  $y(\cdot)\in L^2_{loc}([s-\tau_N,+\infty),\corpsValeurs^d)$ and $u(\cdot)\in L^2_{loc}([s,+\infty),\corpsValeurs^d)$
  satisfy
  \eqref{syst_lin_avec_entree0} for a.e.\ $t \ge s$. Then:

\begin{eqnarray}
\label{eq:variationconstante2}
y_t=U_{2}(t,s)y_s+\displaystyle\int_{s^-}^{+\infty} d_{\alpha}K(t,\alpha)u(\alpha),\qquad  t \ge s.
\end{eqnarray}
\end{proposition}

\begin{pf} 
  Denote by 
  $\eta_t$ the right hand side of
  \eqref{eq:variationconstante2}, and put $\eta(t+\theta)=\eta_t(\theta)$
  for $t\geq s$ with $\theta\in[-\tau_N,0]$.
  Write $\widetilde{y}$ for
  the solution of the homogeneous system \eqref{system_lin_formel} with initial condition $y_s$, so that
$\tilde{y}_t=U_{2}(t,s)y_s$. 
Fix $t\geq s$; for a.e.\ $\theta$ in $[-\tau_N,0]$
satisfying $t+\theta<s$
  (such $\theta$ occur when $t<s+\tau_N$),
  it holds that  $(\int_{s^-}^{+\infty}d_\alpha K(t,\alpha)u(\alpha))(\theta)=\int_{s^-}^{+\infty}d_\alpha X(t+\theta,\alpha)u(\alpha)=0$ (because $X(t+\theta,\alpha)=0$ when $t+\theta<s\leq\alpha$) and that
  $\widetilde{y}_t(\theta)=\widetilde{y}(t+\theta)=y_s(t+\theta)=y(t+\theta)$
  (because
  $\widetilde{y}$ 
  and $y$ 
  have the  same initial condition $y_s$ on $[-\tau_N,0]$).
  Hence, we obtain  indeed that $\eta_t(\theta)=y_t(\theta)$ for a.e.\ $\theta$
  in $[-\tau_N,0]$ such that   $t+\theta<s$.
     Next, for   a.e.\ $\theta \in [-\tau_N,0]$
such that $t+\theta\geq s$,
we deduce from the definition of $\eta_t$,  \eqref{defSOq} (with $q=2$),
\eqref{system_lin_formel} (where $y$ is set to $\widetilde{y}$),
\eqref{defborninf}
and  \eqref{solution_fondamentalev} that 
\begin{align*}
  \eta_t(\theta)&=\tilde{y}(t+\theta)-\int_{s^-}^{+\infty}d_{\alpha}X(t+\theta,\alpha) u(\alpha)\\
  &=\sum\limits_{j=1}^ND_i(t+\theta)\tilde{y}(t+\theta-\tau_j)-
  \sum\limits_{j=1}^ND_i(t+\theta)\int_{s^-}^{+\infty}d_{\alpha}X(t+\theta-\tau_j,\alpha) u(\alpha)+u(t+\theta)\\
&=\sum\limits_{j=1}^ND_j(t+\theta)\eta(t+\theta-\tau_j)+u(t+\theta).
\end{align*}
Altogether, we  proved that to each $t\geq s$ there is a set
$E_t\subset[-\tau_N,0]$ of zero measure such that, for $\theta\in[-\tau_N,0]\setminus E_t$, we have $\eta(t+\theta)=\sum\limits_{j=1}^ND_j(t+\theta)\eta(t+\theta-\tau_j)+u(t+\theta)$ when $t+\theta\geq s$ and $\eta(t+\theta)=y_s(t+\theta)$ when
$t+\theta<s$. Thus, letting $\mathbb{Q}_{+}$ denote the nonnegative rational
numbers, the set $E :=\cup_{q\in\mathbb{Q}_+} \{s+q+E_{s+q}\}$ has
measure zero and since each $t\ge s$ can be written as $s+q+\theta$
with $q\in \mathbb{Q}_{+}$ and $\theta\in[-\tau_N,0]$, we deduce that
$\eta(t)$ satisfies \eqref{syst_lin_avec_entree0}
when $t\geq s$ and $t\notin E$, with initial condition $\eta=y_s$
a.e.\ on $[-\tau_N,0]$. Hence, by uniqueness of a solution with given initial condition, $\eta=y$ a.e.
\qed
\end{pf}

Let us now fix $t$ and consider the function $s\mapsto X(t,t-s)$
from $\corpsR$ to $\corpsValeurs^{d\times d}$. It is right-continuous and lies in
$BV_{loc}(\corpsR,\corpsValeurs^{d\times d})$, moreover it is identically zero for $s<0$.
Thus, on each interval $I\subset\corpsR$ which is bounded on the right,
it generates a $\corpsValeurs^{d\times d}$-valued measure $\nu_{(X(t,t-.)_{|I})}$
which is but the image of
  $d_sX(t,s)$  under the map
  $s\mapsto t-s$ from $t-I$ onto $I$ (where $t-I$ indicates the set of time instants of the form $t-\tau$ for $\tau\in I$).
When $I$ is understood,
we denote this measure by
$d_sX(t,t-s)$ and  we get from \eqref{formmeas} that
\begin{equation}
  \label{formmeasi}
d_sX(t,t-s)=-\sum_{\mathfrak{f}\in\mathcal{F}\cap \check{I}} \mathfrak{C}_\mathfrak{f}(t)\delta_\mathfrak{f},
  \end{equation}
  where $\check{I}$ denotes $I$ deprived from its left endpoint
  (if contained in $I$) and 
  the $\mathfrak{C}_\mathfrak{f}$ are as in \eqref{formmeas}.
  By periodicity of the $D_j$, one sees that $X(t,t-s)$ is periodic in $t$
  and therefore the $\mathfrak{C}_\mathfrak{f}(t)$ are periodic as well.
  If the $D_j$ belong to $C^{1,\delta}(\corpsR,\corpsValeurs^{d\times d})$ then so do
  the   $\mathfrak{C}_\mathfrak{f}$, and we may
   define the measure
\begin{equation}
  \label{formmeasider}
  \frac{\partial}{\partial t}
  d_sX(t,t-s):=-\sum_{\mathfrak{f}\in\mathcal{F}\cap \check{I}} \left(\frac{\partial}{\partial t}\mathfrak{C}_{\mathfrak{f}}(t)\right)\delta_\mathfrak{f}
  \end{equation}
  with coefficients
  in $C^\delta(\corpsR,\corpsValeurs^{d\times d})$.
  The number of terms in the right hand sides of
  \eqref{formmeasi}, \eqref{formmeasider}, and the number of sums of
  products of the $D_j(t-\mathfrak{f}')$ 
involved in these terms, tend to $+\infty$ with the length of $I$.
 Nevertheless, in the next section we show  that the growth of these quantities is at most exponential with that length.

\subsection{A priori estimates}
 For $I=[0,\tau]$, we shall need basic {\it a priori} estimates, independent of $t$,
  for the quantities: 
  \begin{align}
    \label{defn}
    \left\|d_sX(t,t-s)\right\|_{I}:=&\sum_{\mathfrak{f}\in\mathcal{F}\cap\check{I}}\vertiii{\mathfrak{C}_{\mathfrak{f}}(t)},\\
    \left\|\frac{\partial}{\partial t}d_sX(t,t-s)\right\|_{I}:=&\sum_{\mathfrak{f}\in\mathcal{F}\cap\check{I}}\vertiii{\frac{\partial}{\partial t}\mathfrak{C}_{\mathfrak{f}}(t)},
    \label{dfn1}\\
        \Lambda_{I,\delta}\left(\frac{\partial}{\partial t}d_sX(t,t-s)\right):=&
  \sum_{\mathfrak{f}\in\mathcal{F}\cap\check{I}}\Lambda_\delta\left(\frac{\partial}{\partial t}\mathfrak{C}_{\mathfrak{f}}(t)\right),  \label{dfn2}   
    \end{align}
where $\Lambda_\delta(g)$ indicates the H\"older constant of $g\in C^\delta(\corpsR,\corpsValeurs^{d\times d})$.
Note that \eqref{defn} and \eqref{dfn1}  are just the total variations
of the measures $d_sX(t,t-s)$ and $\frac{\partial}{\partial t}d_sX(t,t-s)$
on $I$, respectively.
\begin{proposition}
\label{majoration-sol-fond}
Assume that $D_j \in C^{1,\delta}(\corpsR,\corpsValeurs^{d\times d})$ is $T$-periodic
for $1\leq j\leq N$ and some $\delta \in (0,1)$.
Then, there exist $K,\gamma\geq0$ such that, for all $\tau\geq0$ and $t\in\corpsR$:
\begin{eqnarray}
\label{majoration_sol_fond_equation}
\max\left\{ \left\|d_sX(t,t-s)\right\|_{[0,\tau]},\left\|\frac{\partial}{\partial t}d_sX(t,t-s)\right\|_{[0,\tau]}, \Lambda_{[0,\tau],\delta}\left(\frac{\partial}{\partial t}d_sX(t,t-s)\right)\right\}
  \le K e^{\gamma \tau}.
\end{eqnarray}

\end{proposition}

\begin{pf}

    Observe that $d_sX(t,t-s)$ is the image of the measure
  $d_sX(t,s)$ under the map
  $s\mapsto t-s$ from $[t-\tau,t]$ onto $[0,\tau]$. Hence, it follows from
  \eqref{solution_fondamentalev} that
\begin{align*}
  d_sX(t,t-s)&=
 \sum_{j=1}^N D_j(t) d_sX(t-\tau_j,t-s)\\
  &=\sum\limits_{j=1}^N D_j(t) d_sX(t-\tau_j,t-\tau_j-(s-\tau_j))\qquad \mathrm{on}\  [0,\tau],
\end{align*}
as well as
\begin{align*}
  \frac{\partial}{\partial t}d_sX(t,t-s)
    &=\sum\limits_{j=1}^N \frac{\partial}{\partial t}D_j(t) d_sX(t-\tau_j,t-s)+
  \sum\limits_{j=1}^N D_j(t) \frac{\partial}{\partial t}d_sX(t-\tau_j,t-s)\\
  &\hspace{-2em}=\sum\limits_{j=1}^N \frac{\partial}{\partial t}D_j(t) d_sX(t-\tau_j,t-\tau_j-(s-\tau_j))+
     \sum\limits_{j=1}^N D_j(t) \frac{\partial}{\partial t}d_sX(t-\tau_j,t-\tau_j-(s-\tau_j)).
     \end{align*}
Since $X(t-\tau_j,t-\tau_j-(s-\tau_j))$ is identically zero for $s<\tau_j$,
we deduce from the previous identities that if $K'$ is an upper bound
(made independent of $t$ by periodicity) for the
$\vertiii{D_j}$, the $\vertiii{\frac{\partial}{\partial t}D_j}$,
the $\Lambda_\delta(D_j)$ and the $\Lambda_\delta(\frac{\partial}{\partial t}D_j)$, then
\begin{eqnarray}
 \label{eq:mysteredp}
  \|d_sX(t,t-s)\|_{[0,\tau]}&\!\!\!\!\!\!\!\!\!\!\!\!\!\!\!\!\!\!\!\!\!\!\!
                              \!\!\!\!\!\!\!\!\!\!\!\!\!\!\!\!\!\!\!\!\!\!\!
                              \!\!\!\!\!\!\!\!\!\!\!\!\!\!\!\!\!\!\!\!\!\!\!
                              \!\!\!\!\!\!\!\!\!\!\!\!\!\!\!\!\!\!\!\leq
  K'\sum\limits_{j=1}^N \|d_\beta X(t-\tau_j,t-\tau_j-\beta)\|_{[0,\tau-\tau_j]},\\
  \label{eq:mysteredd}
\!\!\!\!\!\!\!\!\!\!\!\!\!\!\!\!\!  \|\frac{\partial}{\partial t}d_sX(t,t-s)\|_{[0,\tau]}&\!\!\!\!\!\!\!\!\!\leq \!K'\!\sum\limits_{j=1}^N \!\!\|d_\beta X(t-\tau_j,t-\tau_j-\beta)\!\|_{[0,\tau-\tau_j]}\!\!+
                                                         \!\!\|\frac{\partial}{\partial t}d_\beta X(t-\tau_j,t-\tau_j-\beta)\!\|_{[0,\tau-\tau_j]},\\
  \nonumber
  \Lambda_{[0,\tau],\delta}(
  \frac{\partial}{\partial t}d_sX(t,t-s))&\leq \!K'\!\sum\limits_{j=1}^N \Lambda_{[0,\tau-\tau_j],\delta}(X(t-\tau_j,t-\tau_j-\beta))+\|d_\beta X(t-\tau_j,t-\tau_j-\beta)\!\|_{[0,\tau-\tau_j]}\\
  &\!\!\!\!\!\!\!\!\!\!\!\!\!\!\!\!+\Lambda_{[0,\tau-\tau_j],\delta}(\frac{\partial}{\partial t}X(t-\tau_j,t-\tau_j-\beta)) + \|\frac{\partial}{\partial t}d_\beta X(t-\tau_j,t-\tau_j-\beta)\!\|_{[0,\tau-\tau_j]},
    \label{eq:mysteredh}
\end{eqnarray}
where we used that $\Lambda_\delta(AB)\leq\Lambda_\delta(A)\vertiii{B}+
\Lambda_\delta(B)\vertiii{A}$.
The proof may now be completed by an inductive step, similar to
the one used to establish  \eqref{bound_sol_fund}. More precisely,
pick $\gamma$,  $K''$ positive  large enough that 
$K'Ne^{-\gamma \tau_1}<1/4$ and
\begin{equation}
  \label{sec}
  \max\{\vertiii{\mathfrak{C}_{\tau_1}(t)},
  \vertiii{\frac{\partial}{\partial t} \mathfrak{C}_{\tau_1}(t)},\Lambda_\delta(\mathfrak{C}_{\tau_1}),\Lambda_\delta(\frac{\partial}{\partial t}\mathfrak{C}_{\tau_1})\}\leq K'', \qquad t\in\corpsR.
  \end{equation}
Note that such a $K''$ indeed exists, as $\mathfrak{C}_{\tau_1}$
is periodic and lies in $C^{1,\delta}(\corpsR)$. From \eqref{formmeasi}  one sees that $\|d_sX(t,t-s)\|_{[0,\tau]}$
is equal to $0$ when $\tau\in[0,\tau_1)$ and to $\vertiii{\mathfrak{C}_{\tau_1}(t)}$ when $\tau=\tau_1$. {\it A fortiori} then, for $\tau\in[0,\tau_1]$ we have:
\[\max\left\{\|d_sX(t,t-s)\|_{[0,\tau]}, \|\frac{\partial}{\partial t}d_sX(t,t-s)\|_{[0,\tau]}\right\}\leq K'' e^{\gamma\tau},\quad
  \Lambda_{[0,\tau],\delta}(\frac{\partial}{\partial t}d_sX(t,t-s))\leq
  2K'K''e^{\gamma\tau},\]
 so that  \eqref{majoration_sol_fond_equation} holds with $K:=\max\{K'',2K'K''\}$ for $0\leq\tau\leq\tau_1$. Now, assume that \eqref{majoration_sol_fond_equation} is true for all $t$ and 
  $\tau\in [0,k\tau_1)$, with $k>0$ an integer. For $\tau\in(k\tau_1,(k+1)\tau_1]$, we get in view of this hypothesis
  and the definition of $\gamma$, that
  \[\|d_sX(t,t-s)\|_{[0,\tau]}\leq K'NK e^{\gamma{(\tau-\tau_1)}}\leq Ke^{\gamma\tau}/4    \qquad \mathrm{by}\  \eqref{eq:mysteredp},\]
  \[\|\frac{\partial}{\partial t}d_sX(t,t-s)\|_{[0,\tau]}\leq 2K'NK e^{\gamma{(\tau-\tau_1)}}\leq Ke^{\gamma\tau}/2    \qquad \mathrm{by}\  \eqref{eq:mysteredd},\]
  \[\Lambda_{[0,\tau],\delta}(\frac{\partial}{\partial t}d_sX(t,t-s))\leq
    4K'NK e^{\gamma{(\tau-\tau_1)}}\leq Ke^{\gamma\tau}    \qquad \mathrm{by}\  \eqref{eq:mysteredh}.\]

By induction on $k$, this completes the proof. 
\qed
\end{pf}

\section{Proof of Theorem \ref{theorem_hale_generalization}}
\label{proofGHeHa}

\subsection{Sufficiency 
}
\label{sec:sufficiency}

The proof of sufficiency can be  modeled after the one in \cite{Henry1974} for the time-invariant case,
but new technicalities arise because applying Laplace transformation
to linear relations with time-varying
coefficients
 yields functional rather than algebraic equations.
Assuming that both \ref{assumption1}, \ref{assumption2}
hold in
Theorem~\ref{theorem_hale_generalization}
we shall prove exponential stability of System \eqref{system_lin_formel}
in three steps.

\subsection*{Step 1}
By assumption the $D_j$ have H\"older continuous derivatives, so
their Fourier coefficients \eqref{fourier_coeef} satisfy  
\begin{eqnarray}
\label{regularity_fourier_coefficient}
\displaystyle\vertiii{\check{D}_j(k)} \le \frac{C}{1+|k|^{1+\delta}} \mbox{, $j \in \{1,\cdots,n \}$},
\end{eqnarray}
where $C$ is a positive constant and $\delta \in ]0,1[$ is the H\"older exponent of the derivative; indeed, \eqref{regularity_fourier_coefficient} follows at
once from the fact that the modulus of the $k$\textsuperscript{th} Fourier coefficient of $\frac{d}{dt}D_j$
is bounded by $C'/(1+n)^\delta$,
see \cite[Ch.\ 2, thm.\ 4.7]{zygmund2002trigonometric}. Note, for later use,
that the constant $C$ in \eqref{regularity_fourier_coefficient} depends affinely
on the H\"older constant of $\frac{d}{dt}D_j$.
Substituting the Fourier expansion of $D_j$
in
\eqref{solution_fondamentale} yields:
\begin{eqnarray}
\label{eq:18.99}
   X(t,s)=I_d+ \sum\limits_{k \in \mathbb{Z}}    
\sum\limits_{j=1}^N  e^{i k \omega t}\check{D}_j(k) X(t-\tau_j,s) \mbox{ \ if $t \ge s$},\qquad \mbox{$X(t,s)=0$ if $t<s$,} 
\end{eqnarray}
where the right hand side of \eqref{eq:18.99}
is absolutely convergent, locally uniformly in $(t,s)$, thanks to
\eqref{regularity_fourier_coefficient} and the local boundedness of $X$.
Hereafter, we shall denote with a hat the Laplace transform with respect to
the first variable; {\it e.g.}, we put
\begin{eqnarray}
\hat{X}(p,s):=\int_{-\infty}^{+ \infty} e^{-pt} X(t,s)dt.
\end{eqnarray}
Taking the Laplace transform of both sides  of \eqref{eq:18.99}
and interchanging the series and integral signs
(this is possible thanks to
\eqref{regularity_fourier_coefficient}), we obtain 
that for $p\in\corpsC$ with $\Re(p)>\lambda$  where 
$\lambda$ is as in \eqref{bound_sol_fund}:
\begin{equation}
\label{eq:19}
 \begin{split}
    \int_{-\infty}^{+ \infty} e^{-pt} X(t,s)dt=\int_{-\infty}^{+ \infty}  \mathds{1}_{[s,+\infty)}(t) e^{-pt} I_d \,dt + \sum\limits_{k \in \mathbb{Z}}    
\sum\limits_{j=1}^N    \int_{-\infty}^{+ \infty} e^{(-p+i k \omega)t}\check{D}_j(k)X(t-\tau_j,s)  dt , 
    \end{split}
\end{equation}
where $\mathds{1}_{[s,+\infty)}$ is the characteristic function equal to $1$ on $[s,+\infty)$ and to $0$ elsewhere; note that the factor $X(t-\tau_j,s)$ in
the  integrand on the right hand side of \eqref{eq:19} is zero for
$t<s$. 
 By an elementary change of variable,
 we deduce from \eqref{eq:19} that
\begin{eqnarray}
\label{eq:21}
\hat{X}(p,s)=\frac{e^{-ps}}{p}I_d+\sum\limits_{k \in \mathbb{Z}} \sum\limits_{j=1}^N    e^{(-p+i k \omega)\tau_j}\check{D}_j(k) \hat{X}(p-i k \omega,s), \qquad \Re(p)>\lambda.
\end{eqnarray}
This can be rephrased as
\begin{equation}
  \label{eq:7}
  \sum_{k\in\mathbb{Z}}\left(R(p)\right)_{0,k}\hat{X}(p-i k \omega,s)
  =\frac{e^{-ps}}{p}I_d,
\end{equation}
with $R(p)$ 
as in \eqref{defHc}, \eqref{defHc-bis} and where
  $(R(p))_{0,k}$ stands for the block at the intersection of block-line
  $0$ and block-column $k$. Substituting
$p+i n \omega$  for $p$ in \eqref{eq:7} and
letting $n$ range over $\mathbb{Z}$,
we obtain a system of countably many
equations that  may be written as
\begin{eqnarray}
\label{syst_ref}
R(p) \,\hat{\mathbb{X}}(p,s)=\hat{e}(p,s)\,,
\end{eqnarray}
where $\hat{\mathbb{X}}(p,s)$ and $\hat{e}(p,s)$ are the infinite 
$d\times d$ block vectors:
\begin{equation}
  \label{eq:9}
  \hat{\mathbb{X}}(p,s):=\left( \begin{array}{c}
\vdots \\
\hat{X}(p+i \omega,s)\\
\hat{X}(p,s)\\
\hat{X}(p-i  \omega,s) \\
\vdots 
                                \end{array} \right)
\ \ \text{and}\ \ 
\hat{e}(p,s):=\left( \begin{array}{c}
\vdots \\
\frac{e^{-(p+i \omega)s}}{p+i \omega}I_d\\
\frac{e^{-ps}}{p}I_d\\
\frac{e^{-(p-i  \omega)s}}{p-i  \omega} I_d\\
\vdots 
\end{array} \right)\,.
\end{equation}
Clearly, these are the Laplace transforms of  infinite $d\times d$ block vectors $ \mathbb{X}(t,s)$ and $e(t,s)$:
\begin{align}
  \label{eq:10}
  &\hat{\mathbb{X}}(p,s)=\int_{-\infty}^{+ \infty} e^{-pt} \mathbb{X}(t,s) \mathrm{d}t
  \,,\ \
  \hat{e}(p,s) (p,s)=\int_{-\infty}^{+ \infty} e^{-pt} e(t,s) \mathrm{d}t\,,
\\
  \label{eq:8}
  &\text{with}\quad
  \mathbb{X}(t,s):=
     \left( \begin{array}{c}
                            \vdots \\
                            e^{-i \omega t}X(t,s)\\
                            X(t,s)\\
                            e^{i \omega t}X(t,s) \\
                            \vdots 
     \end{array} \right)
\,,
   e(t,s):=
   \left( \begin{array}{c}
            \vdots \\
            I_de^{-i \omega t}\mathds{1}_{[s,+\infty)}\\
            I_d\mathds{1}_{[s,+\infty)}\\
            I_de^{i \omega t}\mathds{1}_{[s,+\infty)}\\
            \vdots 
    \end{array} \right)\,.
\end{align}
The infinite dimensional  linear constant system \eqref{syst_ref}
recasts the finite-dimensional periodic time-varying delay system
\eqref{solution_fondamentale} in terms of Fourier series and Laplace transforms.
In order to estimate $\mathbb{X}(p,s)$ --and eventually $X(t,s)$-- we
shall invert the operator $R(p)$; this is our next step.

\subsection*{Step 2}
For $A:=(a_{i,j})_{i,j \in \mathbb{Z}}$  a doubly infinite block matrix,
where $a_{i,j}$ is a $d \times d$ complex matrix for each $i$ and $j$ in
$\mathbb{Z}$, we say that $A$ has off diagonal decay of order $r$ if there is a constant $C$ such that $\displaystyle\vertiii{a_{i,j}}\leq C(1+|i-j|)^{-r}$.
We also define the Wiener norm of $A$ to be
\begin{eqnarray}
\label{def_norm_Wiener}
\displaystyle\vertiii{A}_{\mathfrak{W}}:=\sum\limits_{k\in \mathbb{Z}} \underset{|i-j| = k}{\sup} \displaystyle\vertiii{a_{i,j}}.
\end{eqnarray}
Note that
 $\displaystyle\vertiii{A}_2\leq \displaystyle\vertiii{A}_{\mathfrak{W}}$,
as follows immediately from the Schur test \cite{brown}.
We denote by $B(\ell^2(\mathbb{Z},\corpsValeurs^d),\ell^2(\mathbb{Z},\corpsValeurs^d))$ the space of
doubly infinite block matrices $A$ such that $\displaystyle\vertiii{A}_2<\infty$,
and by $\mathfrak{W}(\ell^2(\mathbb{Z},\corpsValeurs^d),\ell^2(\mathbb{Z},\corpsValeurs^d))$ the subspace of
those $A$ satisfying $\displaystyle\vertiii{A}_{\mathfrak{W}}<\infty$.
It is easy to check that $\mathfrak{W}(\ell^2(\mathbb{Z},\corpsValeurs^d),\ell^2(\mathbb{Z},\corpsValeurs^d))$  is a subalgebra of $B(\ell^2(\mathbb{Z},\corpsValeurs^d),\ell^2(\mathbb{Z},\corpsValeurs^d))$.

One can see from \eqref{regularity_fourier_coefficient}  and \eqref{defHc-bis}
that
$R(p)$ possesses off-diagonal decay of order $1+\delta$, moreover the constant is uniform over any half-space $\Re(p)\geq a$. From this,
it follows at once that 
$\displaystyle\vertiii{R(p)}_{\mathfrak{W}}$ is uniformly bounded over such a
half-space. Hence, $R(p) \in  \mathfrak{W}(\ell^2(\mathbb{Z},\corpsValeurs^d),\ell^2(\mathbb{Z},\corpsValeurs^d))$ for all $p$ and clearly, $R : \corpsC \longrightarrow B(\ell^2(\mathbb{Z},\corpsValeurs^d),\ell^2(\mathbb{Z},\corpsValeurs^d))$  is a Banach valued holomorphic function.

Now, assumption \ref{assumption1} of Theorem \ref{theorem_hale_generalization} tells us that
$R(p)$ is invertible in  $B(\ell^2(\mathbb{Z},\corpsValeurs^d),\ell^2(\mathbb{Z},\corpsValeurs^d))$
for $\Re(p) \ge \beta$, and assumption \ref{assumption2} that
the inverse operator $R(p)^{-1}$ has uniformly bounded
$\displaystyle\vertiii{\cdot}_2$-norm. 
Therefore, as $R(p)$ has off-diagonal decay of order $1+\delta$,
we get from
\cite[thm.\ 1.2]{grochenig2014norm} and the boundedness of
$\displaystyle\vertiii{R(p)^{-1}}_2$ that $R(p)^{-1}$ also has
off-diagonal decay of order $1+\delta$, uniformly for
 $\Re(p) \ge \beta$. In particular, 
$\displaystyle\vertiii{R(p)^{-1}}_{\mathfrak{W}}$ is uniformly bounded for $ \Re(p)  \ge \beta$:
\begin{eqnarray}
  \label{binvW}
  \mbox{$\displaystyle\vertiii{\left[I_{\infty}- \sum\limits_{i=j}^N e^{- p \tau_j} L_{D_j} \Delta_{\tau_j,\omega}\right]^{-1}}_{\mathfrak{W}}\leq C_1$,}\qquad
p \in \{ z \in \corpsC|\Re(z)\ge\beta\}.\end{eqnarray}

By \eqref{binvW},
applying  to \eqref{defHc} the result of  \cite[thm.\ 2.6]{krishtal2011wiener}
that generalizes the Wiener lemma to Banach algebra-valued almost
periodic functions (the Banach algebra here is $\mathfrak{W}(\ell^2(\mathbb{Z},\corpsValeurs^d),\ell^2(\mathbb{Z},\corpsValeurs^d))$, see 
\cite{BALAN2010339} for a still more general version of that theorem),
we get for any $\alpha$ with $\beta<\alpha<0$ that $R(p)^{-1}$ admits a generalized Fourier expansion on the vertical line $p=\alpha+i \mathbb{R}$ of the form : 
\begin{eqnarray}
\label{foureirserie0}
R(\alpha+i\widetilde{\omega})^{-1}=\sum\limits_{k \in \mathbb{Z}}\tilde{R}^{\{k,\alpha\}} e^{i\beta_k \widetilde{\omega} },\qquad \widetilde{\omega} \in \mathbb{R},
\end{eqnarray}
where $\beta_k$ is real and $\tilde{R}^{\{k,\alpha\}} \in \mathfrak{W}(\ell^2(\mathbb{Z},\corpsValeurs^d),\ell^2(\mathbb{Z},\corpsValeurs^d))$ with 
\begin{eqnarray}
\label{somation0}
\sum\limits_{k \in \mathbb{Z}}\displaystyle\vertiii{\tilde{R}^{\{k,\alpha\}}}_{\mathfrak{W}} <+ \infty.
\end{eqnarray}
Moreover,  following the proof of the theorem in section  3 on Dirichlet's series page 147 of \cite{besicovitch1954almost}, which is possible because Cauchy's theorem holds for Banach valued holomorphic functions (one can apply all the arguments there in weak form by evaluating on a fixed vector in $\ell^2(\mathbb{Z},\corpsValeurs^d)$), we deduce that $\tilde{R}^{\{k,\alpha\}}=R^{\{k\}} e^{\alpha \beta_k}$, where $R^{\{k\}} \in B(\ell^2(\mathbb{Z},\corpsValeurs^d),\ell^2(\mathbb{Z},\corpsValeurs^d))$  is independent of $\alpha$. Thus,
one can rewrite  \eqref{foureirserie0} and \eqref{somation0} as
\begin{eqnarray}
\label{foureirserie}
R(\alpha+i\widetilde{\omega})^{-1}=\sum\limits_{k \in \mathbb{Z}}R^{\{k\}}e^{\alpha \beta_k} e^{\beta_k \widetilde{\omega} i}, \qquad \widetilde{\omega} \in \mathbb{R},\quad \alpha>\beta,
\end{eqnarray}
with
\begin{eqnarray}
\label{somation}
\sum\limits_{k \in \mathbb{Z}}\displaystyle\vertiii{R^{\{k\}}}_{\mathfrak{W}}e^{\alpha \beta_k} <+ \infty.
\end{eqnarray}

In addition, observing from the Neumann series that for $\Re(p)$ large enough one has
\begin{eqnarray}
  \label{VN}
R(p)^{-1}=\sum\limits_{n=0}^{+\infty} \left( \sum\limits_{j=1}^N L_{D_j} \Delta_{\tau_j,\omega} e^{- p \tau_j} \right)^n,
\end{eqnarray}
we deduce that
\begin{eqnarray}
\label{form_beta_k}
\{ \beta_k|k \in \mathbb{Z} \} \subset \{ \sum\limits_{j=1}^N n_j \tau_j| \mbox{ $n_j$ non-positive integers for $j=1,\cdots,N$} \}=-\mathcal{F}
\end{eqnarray}
where $\mathcal{F}$ was defined in \eqref{set_discontinuities}.
Note that \eqref{foureirserie}, \eqref{somation} and \eqref{form_beta_k} are
reminiscent of \cite[p.\ 429, equations (12.15.12) and (12.15.13)]{Bellman1959}, that deals  with complex-valued functions.
We also record for later use the following consequence of \eqref{VN}:
\begin{equation}
  \label{VNeu}
  \vertiii{R(p)^{-1}}_2\leq \frac{1}{1-e^{\Re(\tau_1)}\sum_{j=1}^N\vertiii{L_{D_j}}_2\vertiii{\Delta_{\tau_j,\omega}}_2}  \quad
  \mathrm{for}\quad \Re(p)>\gamma_1:=\frac{\log(\sum_{j=1}^N\vertiii{L_{D_j}}_2\vertiii{\Delta_{\tau_j,\omega}}_2)}{\tau_1}.
  \end{equation}

\subsection*{Step 3}
  We now compute $\mathbb{X}(t,s)$. For this, we apply the Laplace inversion formula (see for example \cite[Ch. 1, Lem. 5.2]{Hale}).  It gives us for
  $\Re(p)=c>\lambda$,
  with $\lambda$ as in \eqref{bound_sol_fund}, that
\begin{eqnarray}
\label{eq_inverse_laplace}
\mathbb{X}(t,s)&=&\lim\limits_{\widetilde{\omega} \rightarrow +\infty} \frac{1}{2\pi i} \int_{c-i\widetilde{\omega}}^{c+i\widetilde{\omega}}R(p)^{-1}\hat{e}(p,s) e^{pt} dp,
\end{eqnarray}
  where the convergence in \eqref{eq_inverse_laplace} is understood component-wise; that is to say, for all indices $j \in \mathbb{Z}$, we have:
\begin{eqnarray}
  \lim\limits_{\widetilde{\omega} \rightarrow +\infty}  \vertiii{ \left[ \mathbb{X}(t,s)-\frac{1}{2\pi i} \int_{c-i\widetilde{\omega}}^{c+i\widetilde{\omega}}R(p)^{-1}\hat{e}(p,s) e^{pt} dp  \right]_j}=0.
\end{eqnarray}


We shall need two lemmas, the proof of which is postponed until the end of this section. The first one goes as follows.

\begin{lem}
\label{lemma1}
Let $\beta<\alpha<0$ and $\rho_m:=\frac{2\pi m}{T}+\frac{\pi }{T}$ for $m$
a positive integer. Then, it holds  that
 \begin{eqnarray}
 \label{eq999.99}
 \mathbb{X}(t,s)&=& \lim\limits_{m \rightarrow +\infty} \frac{1}{2\pi i} \int_{\alpha-i\rho_m}^{\alpha+i\rho_m} R(p)^{-1}\hat{e}(p,s) e^{pt} dp+Q(t),
 \end{eqnarray}
where 
\begin{eqnarray}
  \label{defQu}
Q(t):= \sum_{k\in \mathbb{Z}} R(i \omega k)^{-1}_k e ^{i \omega k},
\end{eqnarray}
with $R(p)_k^{-1}$ to mean the $k$\textsuperscript{th} column of $R(p)^{-1}$ for $p\in \mathbb{C}$.
\end{lem}

Assume Lemma \ref{lemma1} for a while and recall that $X(t,s)$ is the  element of index $0$ of
the infinite column vector $\mathbb{X}(t,s)$. Thus,
if for each $k \in \mathbb{Z}$ we let $(R_{0,n}^{\{k\}})_{n \in \mathbb{Z}}$
enumerate the line of index $0$ of the matrix $R^{\{k\}}$ defined in  \eqref{foureirserie}, it is a consequence of \eqref{eq999.99} that
\begin{eqnarray}
\label{eq1000}
X(t,s)&=& \frac{1}{2\pi i} \lim\limits_{m \rightarrow +\infty} \int_{\alpha-i\rho_m}^{\alpha+i\rho_m}\sum\limits_{k,n \in \mathbb{Z}}R_{0,n}^{\{k\}} \frac{e^{-(p-in \omega)s}}{p-i n \omega}e^{p(t+\beta_k)}dp +Q(t)_0,
\end{eqnarray}
where $Q(t)_0$ stands for the element of index $0$ of the column vector $Q(t)$ and the convergence is with respect to the norm $\vertiii{\cdot}$;  note that the convergence of the series in \eqref{eq1000} ensues from
\eqref{somation}.
The second lemma we need is:

\begin{lem}
\label{lemma2}
For all $t$ and $s$ such that $t+\beta_k-s \neq 0$ for all $k \in \mathbb{Z}$,
we have that 
\begin{eqnarray}
\label{eq1003}
X(t,s)=\sum\limits_{(t+\beta_k-s)<0, \,n \in \mathbb{Z}}R_{0,n}^{\{k\}}  e^{i n \omega(t+\beta_k)}+Q(t)_0.
\end{eqnarray}
\end{lem}
Assume Lemma \ref{lemma2} as well for the moment being,
 and note that on a bounded interval can lie only finitely many $\beta_k$,
because the integer coefficients in \eqref{form_beta_k} are non-positive.
Observe also from \eqref{solution_fondamentale} that $X(t,.)$ is left
continuous. Hence, Lemma \ref{lemma2} allows us to compute
$X(t,s)$ for all $(t,s)$. Hereafter,
we fix the symbol $s$ to mean the
initial time,  and  we evaluate the variation of
$X(t,\tau)$ in its second argument when the latter
ranges over $[s,s+\tau_N]$. To this effect, we substitute $\tau$ for $s$ in
\eqref{eq1003} and observe,
since $Q(t)_0$ is independent of $\tau$ while the sum in \eqref{eq1003}
is piecewise constant in $\tau$ with jumps (induced by the constraint on indices) at the $t+\beta_k$
 where it is left continuous, that
\begin{eqnarray}
\label{majoration_variation}
W_{[s,s+\tau_N]}X (t,.) &\le& \sum\limits_{ 0 \le t+\beta_k-s \le \tau_N, n \in \mathbb{Z}}\displaystyle\vertiii{R_{0,n}^{\{k\}}} \\
                                                       \label{majoration_variation1} &\le& \left( \sum\limits_{k,n \in \mathbb{Z}}\displaystyle\vertiii{R_{0,n}^{\{k\}}}
                                                              e^{\alpha \beta_k}e^{|\alpha|\tau_N} \right) e^{\alpha(t-s)},
\end{eqnarray}
where we used in the second inequality that $\alpha(t+\beta_k-s-\tau_N)\ge 0$
because $\alpha<0$. If $K'>0$ is an upper bound for
  $\{\vertiii{D_i(t)}, \, 1\leq i\leq N,\,t\in\corpsR\}$, then  we deduce
from \eqref{formula_representation} and \eqref{majoration_variation1}
that
\begin{eqnarray}
  \|y(t)\| &\le& K'N\left(
                 W_{[s,s+\tau_N]}X (t,.)\right) \|\phi\|_{C^0} \\
&\le& K'N  \left( \sum\limits_{k,n \in \mathbb{Z}}\displaystyle\vertiii{R_{0,n}^{\{k\}}}e^{\alpha \beta_k}e^{|\alpha|\tau_N} \right) e^{\alpha(t-s)} \|\phi\|_{C^0},
\end{eqnarray}
thereby showing that System \eqref{system_lin_formel} is $C^0$-exponentially stable.
This achieves the proof of sufficiency in
Theorem~\ref{theorem_hale_generalization}
granted Lemmas \ref{lemma1} and \ref{lemma2}, that we now establish.


\begin{pflemma1}
Considering the subsequence $\rho_m=\frac{2 \pi m}{T}+\frac{\pi}{T}$, we deduce from \eqref{eq_inverse_laplace} that :
\begin{eqnarray}
\label{eq998.99}
\mathbb{X}(t,s)&=&\lim\limits_{m \rightarrow +\infty} \frac{1}{2\pi i} \int_{c-i\rho_m}^{c+i\rho_m}R(p)^{-1}\hat{e}(p,s) e^{pt} dp.
\end{eqnarray}
As each component of $R(p)^{-1}\hat{e}(p,s)$ is a meromorphic function in the half plane $\{p \in \corpsC| \Re(p) >\beta \}$
with simple poles in the set $\{{i k\pi}/{T}:\,k \in \mathbb{Z} \}$, by
\eqref{eq:9} and the assumed conditions $(i)$,  $(ii)$
of Theorem~\ref{theorem_hale_generalization},
we get from  \eqref{eq998.99} together with  the residue theorem that
\begin{eqnarray}
\label{eq999}
\mathbb{X}(t,s)=\lim\limits_{m \rightarrow +\infty} \frac{1}{2\pi i}  \left( \int_{\alpha-i\rho_m}^{\alpha+i\rho_m}+\int_{\alpha+i\rho_m}^{c+i\rho_m} -\int_{\alpha-i\rho_m}^{c-i\rho_m}  \right) R(p)^{-1}\hat{e}(p,s) e^{pt}dp +Q(t)
\end{eqnarray}
where $Q$ is as in \eqref{defQu};
observe that taking  the limit in \eqref{eq999} is indeed permitted, since
\begin{eqnarray}
\|Q(t)\|_2 \le \displaystyle\vertiii{R(0)^{-1}}_{\mathfrak{W}}<+\infty.
\end{eqnarray}
To establish the lemma, it remains to show  that 
\begin{eqnarray}
\label{eq95657}
\lim\limits_{m \rightarrow +\infty}\int_{\alpha+i\rho_m}^{c+i\rho_m}R(p)^{-1}\hat{e}(p,s) e^{pt}dp=\lim\limits_{m \rightarrow +\infty}\int_{\alpha-i\rho_m}^{c-i\rho_m} R(p)^{-1}\hat{e}(p,s) e^{pt}dp =0.
\end{eqnarray}
For this, let us write the entry with index $j\in\ZZ$ of the block vector $\int_{\alpha+i\rho_m}^{c+i\rho_m}R(p)^{-1}\hat{e}(p,s) e^{pt}dp$
as
\begin{eqnarray}
  \label{repind}
\int_{\alpha+i\rho_m}^{c+i\rho_m}\sum\limits_{n \in \mathbb{Z}}R(p)^{-1}_{j,n} \frac{e^{-(p-i n \omega)s}}{p-i n \omega} e^{pt} dp, 
\end{eqnarray}
where $R(p)^{-1}_{j,n}$ denotes the block entry with index $(j,n)$ of $R(p)^{-1}$. As pointed out in step 2,  $R(p)^{-1}$ has
off-diagonal decay of order $1+\delta$, uniformly for
$\Re(p) \ge \beta$. Thus, if we pick $\varepsilon>0$,
there is $\mathfrak{n}=\mathfrak{n}(j)$ such that
$\sum\limits_{|n|\geq\mathfrak{n}}\vertiii{R(p)^{-1}_{j,n}}<\varepsilon$,
uniformly with respect to $\Re(p)\in[\alpha,c]$, and since 
\begin{eqnarray}
\label{eq95659}
 \mbox{$|\rho_m -n \omega| \ge \frac{\pi}{T}$ for all $n$ and $m$}
 \end{eqnarray}
 we have that
\begin{equation}
\label{eq956581}
\displaystyle\vertiii{ \int_{\alpha+i\rho_m}^{c+i\rho_m}\sum\limits_{|n|\geq\mathfrak{n}}R(p)^{-1}_{j,n} \frac{e^{-(p-in \omega)s}}{p-in\omega} e^{pt} dp }\le  \varepsilon \,e^{c (t-s)}\,\frac{(c-\alpha)T}{\pi}
\end{equation}
which is arbitrary small with $\varepsilon$. In another connection,
for fixed $j$ and $n\in\mathbb{Z}$, we get since $|p-in \omega|$ becomes arbitrary large with $m$
while the other terms in the integrand
are uniformly bounded with respect to
$\Re(p)\in[\alpha,c]$ that
 \begin{equation}
\label{eq956582}
\underset{m\to\infty}{\lim} \int_{\alpha+i\rho_m}^{c+i\rho_m}R(p)^{-1}_{j,n} \frac{e^{-(p-in \omega)s}}{p-in \omega} e^{pt} dp =0.
\end{equation}

In view of \eqref{eq956581} and \eqref{eq956582}, we conclude
 that 
\begin{eqnarray}
\lim\limits_{m \rightarrow +\infty}\int_{\alpha+i\rho_m}^{c+i\rho_m}R(p)^{-1}\hat{e}(p,s) e^{pt}dp=0,
 \end{eqnarray}
and in the same way one can prove that
 \begin{eqnarray}
 \lim\limits_{m \rightarrow +\infty}\int_{\alpha-i\rho_m}^{c-i\rho_m} R(p)^{-1}\hat{e}(p,s) e^{pt}dp =0,
 \end{eqnarray}
 as wanted.
 \qed
\end{pflemma1}

\begin{pflemma2}
  Let $t$ and $s$ be such that $t+\beta_k-s \neq 0$ for all $k\in\ZZ$. By \eqref{eq95659} and \eqref{somation} that allow us to integrate termwise the series below, we get on performing the change of variable $p \rightarrow p-i n \omega$ in the integrals corresponding to index $n$ that
  \begin{equation}
\label{eq1001.2}
\int_{\alpha-i\rho_m}^{\alpha+i\rho_m} \sum\limits_{k,n \in \mathbb{Z}}R_{0,n}^{\{k\}} \frac{e^{-(p-i n \omega)s}}{p-i n \omega}e^{p(t+\beta_k)}dp =
\sum\limits_{k,n \in \mathbb{Z}}R_{0,n}^{\{k\}} e^{i n \omega(t+\beta_k)}\int_{\alpha-i n \omega-i \rho_m}^{\alpha-i n \omega+i \rho_m} \frac{e^{p(t+\beta_k-s)}}{p}dp.
\end{equation}


Integrating by parts, we \b{obtain:}
\begin{eqnarray}
\label{eq_00.01}
\int_{\alpha-i n \omega-i \rho_m}^{\alpha-i n \omega+i \rho_m} \frac{e^{p(t+\beta_k-s)}}{p}dp=\left[\frac{e^{p(t+\beta_k-s)}}{p(t+\beta_k-s)} \right]^{\alpha-i n \omega+i \rho_m}_{\alpha-i n \omega-i \rho_m}+\int_{\alpha-i n \omega-i \rho_m}^{\alpha-i n \omega+i \rho_m} \frac{e^{p(t+\beta_k-s)}}{(t+\beta_k-s)p^2}dp,
\end{eqnarray}
and by \eqref{form_beta_k} there is  $\delta>0$ such that $|t+\beta_k-s|\ge \delta$ for all $k \in \mathbb{Z}$. Hence,
using \eqref{eq95659}, we obtain:

\begin{eqnarray}
\label{eq_00.1}
\left|\left[\frac{e^{p(t+\beta_k-s)}}{p(t+\beta_k-s)} \right]^{\alpha-i n \omega+i \rho_m}_{\alpha-i n \omega-i \rho_m}\right| &\le& \frac{1}{|t+\beta_k-s|} \left( \frac{e^{\alpha(t+\beta_k-s)}}{\sqrt{\alpha^2+ \left(\rho_m-i n \omega \right)^2}}+\frac{e^{\alpha(t+\beta_k-s)}}{\sqrt{\alpha^2+ \left(\rho_m+i n \omega \right)^2}} \right) \nonumber \\
&\le& \frac{2e^{\alpha (t+\beta_k-s)}}{\delta  \sqrt{\alpha^2+ \left(\frac{\pi}{T} \right)^2}}
\end{eqnarray}
and
\begin{eqnarray}
\label{eq_00.2}
\left|\int_{\alpha-i n \omega-i \rho_m}^{\alpha-i n \omega+i \rho_m} \frac{e^{p(t+\beta_k-s)}}{(t+\beta_k-s)p^2}dp\right| &\le& \frac{e^{\alpha(t+\beta_k-s)}}{|t+\beta_k-s|}\int_{-\rho_m- n \omega}^{\rho_m- n \omega} \frac{1}{\alpha^2+p_2^2}dp_2 \nonumber\\
&\le&  \frac{e^{\alpha(t+\beta_k-s)}}{|t+\beta_k-s|} \left[\frac{1}{\alpha} arctan \left( \frac{p_2}{\alpha} \right) \right]_{-\rho_m- n \omega}^{\rho_m- n \omega} \nonumber\\
&\le& \frac{ \pi}{ |\alpha| \delta} e^{\alpha (t+\beta_k-s)}.
\end{eqnarray}

Altogether, we deduce from \eqref{eq_00.01}, \eqref{eq_00.1} and
\eqref{eq_00.2} that there exists a constant $K>0$, independent of $k$, $n$ and $\rho_m$, such that :
\begin{eqnarray}
\label{eq_00.3}
\left| \int_{\alpha-i n \omega-i \rho_m}^{\alpha-i n \omega+i \rho_m} \frac{e^{p(t+\beta_k-s)}}{p}dp\right| \le Ke^{\alpha(t+\beta_k-s)}.
\end{eqnarray}
Consequently, it holds that
\begin{eqnarray}
\label{eq_00.4}
\displaystyle\vertiii{ 
  R_{0,n}^{\{k\}} e^{i n \omega(t+\beta_k)}\int_{\alpha-i n \omega-i \rho_m}^{\alpha-i n \omega+i \rho_m} \frac{e^{p(t+\beta_k-s)}}{p}dp} \le Ke^{\alpha(t-s)} 
  \displaystyle\vertiii{R_{0,n}^{\{k\}}}e^{\alpha \beta_k},
\end{eqnarray}
and since the right hand side of  \eqref{eq_00.4} is summable over $k,n$
because of \eqref{somation}), the dominated convergence theorem allows
us to take the limit termwise
in the right hand side of \eqref{eq1001.2} as $m\to\infty$:

\begin{eqnarray}
  \label{limtmw}
\lim\limits_{m \rightarrow +\infty} \frac{1}{2i \pi} \int_{\alpha-i\rho_m}^{\alpha+i\rho_m}\sum\limits_{k,n \in \mathbb{Z}}R_{0,n}^{\{k\}}  \frac{e^{-(p-i n \omega)s}}{p-i n \omega}e^{p(t+\beta_k)}dp \nonumber \\
 = \sum\limits_{k,n \in \mathbb{Z}} \frac{1}{2i \pi}R_{0,n}^{\{k\}} e^{i n \omega(t+\beta_k)}\lim\limits_{m \rightarrow +\infty} \int_{\alpha-i n \omega-i\rho_m}^{\alpha-i n \omega +i\rho_m} \frac{e^{p(t+\beta_k-s)}}{p}dp.
\end{eqnarray}

Since $\alpha<0$, we get on the one hand by Cauchy's theorem that for
  $\kappa<\alpha$:
    \begin{eqnarray}
  \label{eq100a}
  \int_{\alpha-i n \omega-i\rho_m}^{\alpha-i n \omega+i\rho_m} \frac{e^{p(t+\beta_k-s)}}{p}dp
  =  \left( \int_{\kappa-in\omega+i\rho_m}^{\alpha-in\omega+i\rho_m}+\int_{\kappa-in\omega-i\rho_m}^{\kappa-in\omega+i\rho_m} +\int_{\alpha-in\omega-i\rho_m}^{\kappa-in\omega-i\rho_m}  \right) \frac{e^{p(t+\beta_k-s)}}{p}dp,
  \end{eqnarray}
   and on the other hand by the residue theorem that for $\kappa>0$:
  \begin{eqnarray}
  \label{eq100d}
  \int_{\alpha-i n \omega-i\rho_m}^{\alpha-i n \omega+i\rho_m} \frac{e^{p(t+\beta_k-s)}}{p}dp
  =  2i\pi+\left( \int_{\kappa-in\omega+i\rho_m}^{\alpha-in\omega+i\rho_m}+\int_{\kappa-in\omega-i\rho_m}^{\kappa-in\omega+i\rho_m} +\int_{\alpha-in\omega-i\rho_m}^{\kappa-in\omega-i\rho_m}  \right) \frac{e^{p(t+\beta_k-s)}}{p}dp.
\end{eqnarray}
Since $|p|$ goes to $\infty$ with $m$
while the integrands are uniformly bounded 
with $\Re(p)$, the first and third integrals in the right hand side of
\eqref{eq100a} and \eqref{eq100d} tend to $0$ as $m\to\infty$ for fixed $\kappa$. Hence,
\begin{eqnarray}
  \label{eq100b}
  \lim\limits_{m \rightarrow +\infty}\int_{\alpha-i n \omega-i\rho_m}^{\alpha-i n \omega +i\rho_m} \frac{e^{p(t+\beta_k-s)}}{p}dp
  =  \lim\limits_{m \rightarrow +\infty}\int_{\kappa-in\omega-i\rho_m}^{\kappa-in\omega+i\rho_m}  \frac{e^{p(t+\beta_k-s)}}{p}dp,\qquad \kappa<\alpha,\\
  \label{eq100bb}
   \lim\limits_{m \rightarrow +\infty}\int_{\alpha-i n \omega-i\rho_m}^{\alpha-i n \omega +i\rho_m} \frac{e^{p(t+\beta_k-s)}}{p}dp
  = 2i\pi+ \lim\limits_{m \rightarrow +\infty}\int_{\kappa-in\omega-i\rho_m}^{\kappa-in\omega+i\rho_m}  \frac{e^{p(t+\beta_k-s)}}{p}dp,\qquad 0<\kappa.
\end{eqnarray}
Integrating by parts yields that
\[
  \int_{\kappa-in\omega-i\rho_m}^{\kappa-in\omega+i\rho_m}  \frac{e^{p(t+\beta_k-s)}}{p}dp
  =\frac{e^{\kappa(t+\beta_k-s)}}{i(t+\beta_k-s)}\left[\frac{e^{i(t+\beta_k-s)y}}{\kappa+iy} \right]_{y=-n\omega-\rho_m}^{y=-n\omega+\rho_m}-
  \int_{\kappa-in\omega-i\rho_m}^{\kappa-in\omega+i\rho_m}  \frac{e^{p(t+\beta_k-s)}}{
    i(t+\beta_k-s)p^2}dp
\]
and therefore, if $(t+\beta_k-s)>0$ (resp.\  $(t+\beta_k-s)<0$), then the right hand side
of \eqref{eq100b} (resp.\ \eqref{eq100bb}) goes to 
$0$ (resp.\ $2i\pi$) when $\kappa\to-\infty$ (resp.\ $+\infty$). Thus,
\begin{eqnarray}
\label{eq1002}
\sum\limits_{k,n \in \mathbb{Z}}\frac{1}{2i \pi} R_{0,n}^{\{k\}} e^{i n \omega(t+\beta_k)}\lim\limits_{m \rightarrow +\infty} \int_{\alpha-in\omega-i\rho_m}^{\alpha-in\omega+i\rho_m} \frac{e^{p(t+\beta_k-s)}}{p}dp=\sum\limits_{t+\beta_k-s<0, n \in \mathbb{Z}}R_{0,n}^{\{k\}} e^{i n \omega(t+\beta_k)}.
\end{eqnarray}
In view of \eqref{eq1000}, \eqref{limtmw}
and \eqref{eq1002}, this achieves the proof of Lemma~\ref{lemma2}.
\qed
\end{pflemma2}

\subsection{Necessity 
}
\label{sec:necessariness}

We now assume exponential stability of System \eqref{system_lin_formel}, and prove that conditions \ref{assumption1} and \ref{assumption2}
of Theorem \ref{theorem_hale_generalization} are satisfied.
In the time-invariant case,
necessity can be proven using the spectral theory of semigroups \cite{Henry1974}. Unfortunately,  there seems to be no adequate
notion of semigroup for time-varying coefficients, and we shall
  follow a different path: since exponential stability of System \eqref{system_lin_formel} implies in particular its $L^2$ exponential stability (see
  Proposition~\ref{equi_L2_C0}),
  we may resort to the spectral theory of
monodromy operators for periodic evolution families \cite{buse2001individual}.
These are just solution operators
over a period; {\it i.e.},  if $U_2(.,.)$ refers to the solution operator~\eqref{defSOq} with $q=2$, then the monodromy operator of system \eqref{system_lin_formel}
is $U_2(T,0)$. Choosing zero as initial time  is  arbitrary, and we might
put $U_2(T+s_0,s_0)$  for the monodromy operator with  arbitrary $s_0$, but
even though different choices of $s_0$ lead to different monodromy operators, the
relation of their spectrum to exponential stability is the same, and we stick to $s_0=0$ for definiteness. The result we need is the following.

\begin{proposition}
\label{stability_monodromy}
System \eqref{system_lin_formel} is $L^2$-exponentially stable if and only if the spectral radius of the monodromy operator $U_2(T,0)$ is strictly less than one.
\end{proposition}

\begin{pf}
  This is
  formally a consequence of
 \cite[lem.\ 4.2]{buse2001individual} if definition
 \eqref{eq:stabLq} of exponential stability (with $q=2$) is restricted to $t\geq s\geq0$. However,
 since the constants in that proof only
depend on the $D_j(t)$ through
 $\sup\{\vertiii{D_i(t)}, \, 1\leq i\leq N,\,t\in\corpsR\}$, the case
 of arbitrary $s$ follows by periodicity.
 \qed
 \end{pf}
 
 We mention that  \cite{buse2001individual}
 deals with much more general periodic evolution families than
 \eqref{system_lin_formel}.



%
%

In view of Proposition \eqref{stability_monodromy}, we may 
assume that the spectral radius of the monodromy operator $U_2(T,0)$ is
strictly less than 1.
Then, we will resort to  a classical construction
from Control Theory (called \emph{realization}) to connect the spectrum
of the monodromy operator with  the singular set of the harmonic transfer function.
Specifically, we shall realize
the periodic input-output system \eqref{syst_lin_avec_entree0} as a discrete-time system in state space form, much like in
the reference  \cite{Mollerstedt2000} except that the state is now infinite dimensional.
Consider \eqref{syst_lin_avec_entree0} with initial data $0$
at time $0$, and let $u \in L^2([0,+\infty),\corpsValeurs^d)$ be
an input generating the output  $y(t) \in \corpsValeurs^d$:
   \begin{eqnarray}  
  \label{syst_lin_avec_entreebis}
     y(t)=\sum\limits_{j=1}^ND_j(t)y(t-\tau_j)+u(t) \mbox{,$\qquad$ a.e.\ $t \ge 0$},
     \quad y(t)=0 \ \mathrm{for}\ t\le0.
   \end{eqnarray}
Using \eqref{eq:variationconstante2} with $s=0$, $y_s=0$, and
performing the change of variables $\beta= t-\alpha$, we get that
\begin{eqnarray}
\label{inp-output_equation}
  y(t) 
  =\int_{-\infty}^{t^+} d_{\beta}X(t,t-\beta)u(t-\beta),\qquad \mbox{a.e.\  $t \ge 0$}. 
\end{eqnarray}
Note, since $X(t,t-\beta)$ is identically zero when $\beta<0$, that the
integral in \eqref{inp-output_equation} can also be written as
$\int_{c^\pm}^{t^+} d_{\beta}X(t,t-\beta)u(t-\beta)$ for any $c<0$.
From \eqref{inp-output_equation} and \eqref{formmeasi}, it follows that
\begin{equation}
  \label{vcr}
y(t)=\sum_{\mathfrak{f}\in\mathcal{F}\cap [0,t]} \mathfrak{C}_{\mathfrak{f}}(t) u(t-\mathfrak{f}),\qquad \mathrm{a.e.}\ t\geq0.
 \end{equation}
 Applying to \eqref{vcr} the Schwarz inequality and using
 \eqref{defn}, \eqref{majoration_sol_fond_equation} together with
 the fact that $\mathfrak{C}_0=-I_d$ (as follows readily from
 \eqref{formmeas} and \eqref{solution_fondamentale}), we obtain:
\begin{eqnarray}
    |y(t)|^2\leq&\!\!\!\!\!\!\!\!\left(\sum_{\mathfrak{f}\in\mathcal{F}\cap [0,t]} \vertiii{\mathfrak{C}_{\mathfrak{f}}(t)}^2\right)\left(\sum_{\mathfrak{f}\in\mathcal{F}\cap [0,t]} \| u(t-\mathfrak{f})\|^2\right)\nonumber\\
  \leq& \left(\sum_{\mathfrak{f}\in\mathcal{F}\cap [0,t]} \vertiii{\mathfrak{C}_{\mathfrak{f}}(t)}\right)^2
        \left(\sum_{\mathfrak{f}\in\mathcal{F}\cap [0,t]} \| u(t-\mathfrak{f})\|^2\right)\nonumber\\
  \leq& (1+Ke^{\gamma t})^2\left(\sum_{\mathfrak{f}\in\mathcal{F}\cap [0,t]} \| u(t-\mathfrak{f})\|^2\right)\label{estexpy},
  \qquad \mathrm{a.e.}\ t\geq0.
\end{eqnarray}
In another connection, the cardinality of $\mathcal{F}\cap[0,t]$
is no larger than the number of
$N$-tuples $(q_1,\ldots,q_N)\in\mathbb{N}^N $ satisfying $\sum_{i=1}^N q_i\leq t/\tau_1$ (recall $\tau_1$ is the smallest delay), and so
\begin{equation}
  \label{eq:249a}
  \mathrm{Card}\,\{\mathcal{F}\cap[0,t]\}\leq \left(1+\frac{t}{\tau_1}\right)^{N}\,,\qquad\ t\in[0,+\infty),
\end{equation}
where $\mathrm{Card}\{E\}$ stands for the cardinality of $E$.
Integrating \eqref{estexpy} from $0$ to $\tau>0$ and taking
\eqref{eq:249a} into account, we deduce that
\begin{equation}
  \label{majexpL2py}
\|y\|_{L^2([0,\tau])}\leq \tau^{1/2}\left(1+Ke^{\gamma \tau}\right)\left(1+\frac{\tau}{\tau_1}\right)^{N/2}\|u\|_{L^2([0,+\infty])},
  \end{equation}
  implying that to each $\gamma_1>\gamma$ there is $K_1>0$ for which   
  $\|y\|_{L^2([0,\tau]}\leq K_1 e^{\gamma_1 \tau}\|u\|_{L^2([0,+\infty])}$.
Since $\|u\|_{L^2([0,+\infty])}<+\infty$, this warrants the existence, for $\Re(p)>\gamma$,
of the Laplace transforms of $y(t)$ and $u(t)$:
\begin{eqnarray}
\label{eq:Laplace_transform_input_output}
\hat{Y}(p):=\int_{0}^{+ \infty} e^{-pt} y(t)dt,\qquad
\hat{U}(p):=\int_{0}^{+ \infty} e^{-pt} u(t)dt,
\end{eqnarray}
for if $\Re(p)>\gamma_1>\gamma$ then  $\int_{k}^{k+1} e^{-pt} |u|(t)dt$ and
$\int_{k}^{k+1} e^{-pt} |y|(t)dt$ decay exponentially fast with $k\in\mathbb{N}$,
by the Schwarz inequality.
Proceeding on \eqref{syst_lin_avec_entreebis} 
like we did on \eqref{solution_fondamentale} to obtain \eqref{syst_ref},
namely expanding  the $D_j(t)$ into Fourier series 
and taking Laplace transforms termwise using  
\eqref{regularity_fourier_coefficient}, we get on replacing the complex
variable $p$ by by $p+i n \omega$ and letting $n$ range over
$\mathbb{Z}$ that
\begin{eqnarray}
\label{syst_ref0}
R(p)\hat{\mathbb{Y}}(p)=\hat{\mathbb{U}}(p),\qquad \Re(p) >\gamma,
\end{eqnarray}
where $\hat{\mathbb{Y}}(p)$ and $\hat{\mathbb{U}}(p)$ are given by \eqref{eq:6}
and $R(p)$ by \eqref{defHc}.
Moreover, we see from  \eqref{VNeu} that $R(p)$ is continuously
invertible for
$\Re(p)>\gamma_1$,
so we obtain with  $\gamma_2:=\max\{\gamma,\gamma_1\}$ that
\begin{eqnarray}
\label{syst_ref1}
 \hat{\mathbb{Y}}(p)=R(p)^{-1} \hat{\mathbb{U}}(p),\qquad \Re(p) >\gamma_2.
\end{eqnarray}
Equation~\eqref{syst_ref1} identifies with Equation~\eqref{link_fond_HTF},
showing that the operator-valued holomorphic function
$p\mapsto R(p)^{-1}$ is the Harmonic Transfer Function of System~\eqref{syst_lin_avec_entreebis}; see Definition~\ref{def:HTF}.

 Next, we define the instantaneous transfer function:
\begin{eqnarray}
\label{eq_ITF}
G(t,p)=\int_{-\infty}^{+ \infty} -d_{\tau}X(t,t-\tau)e^{-p \tau} , \qquad t\in\corpsR ,\qquad p \in \corpsC ,\qquad \Re(p)> \gamma,
\end{eqnarray}
where  the integral is understood as $-\lim_{b\to+\infty}\int_{c^\pm}^{b^+}d_{\tau}X(t,t-\tau)e^{-p \tau}$ for any $c<0$. The existence of the limit is guaranteed by Proposition~\ref{majoration-sol-fond}, and the fact that it does not depend
on $c<0$ can be  argued similarly
to the independence of \eqref{defborninf} from the exact value of $b$.
Observe from Proposition~\ref{majoration-sol-fond}
that $t \mapsto  G(t,p)$ lies in $C^{1+\delta}(\corpsR,\corpsValeurs^{d\times d})$, and obviously 
it is  $T$-periodic. Thus, we may expand this function
in Fourier series as
\begin{eqnarray}
\label{fourier_serie_ITF}
G(t,p)= \sum\limits_{k \in \mathbb{Z}}G_k(p) e^{i \omega k t},\qquad \Re(p)>\gamma,
\end{eqnarray} 
where we recall that $\omega:=2 \pi/T$, and the series converges absolutely
for fixed $p$ by a standard estimate used already in
\eqref{regularity_fourier_coefficient}; more precisely,
\cite[Ch.\ 2, thm.\ 4.7]{zygmund2002trigonometric} implies that
\begin{eqnarray}
\label{prop_majoration_coeff_HTF}
\vertiii{G_k(p)} \le \frac{\mathfrak{K}(\Re(p))}{1+|k|^{1+ \delta}},
\end{eqnarray}
where $\mathfrak{K}(\Re(p))>0$ depends on $\Re(p)$, and
also on $K$, $\gamma$ in \eqref{majoration_sol_fond_equation} though we do
not show the latter dependence. Now, for $\mathfrak{f}\in\mathcal{F}$ and $k\in\mathbb{Z}$,
let  $c_{k,\mathfrak{f}}$ be the $k$\textsuperscript{th} Fourier coefficient of the
$T$-periodic function $t\mapsto\mathfrak{C}_{\mathfrak{f}}(t)$.
Since $\mathfrak{C}_{\mathfrak{f}}\in C^{1,\delta}(\corpsR)$, we have the estimate
\begin{eqnarray}
\label{prop_majoration_coeff_HTFc}
\vertiii{c_{k,\mathfrak{f}}} \le \frac{K_{\mathfrak{f}}}{1+|k|^{1+ \delta}}
\end{eqnarray}
so that
$\mathfrak{C}_{\mathfrak{f}}(t)=\sum_{k\in\mathbb{Z}}c_{k,\mathfrak{f}}e^{i\omega kt}$, where the series
is uniformly absolutely convergent. Summing up over those $\mathfrak{f}\in\mathcal{F}$
such that $\mathfrak{f}\leq \tau$, we deduce from \eqref{formmeasi}
that for any $\varepsilon>0$:
\begin{eqnarray}
\label{fact1}
  d_sX(t,t-s)= \sum\limits_{k \in \mathbb{Z}}d\mu_{k,\tau} e^{i \omega k t}\qquad
  \mathrm{on}\quad [-\varepsilon,\tau],
\end{eqnarray}
where the measure $\mu_{k,\tau}$ is equal to
$-\sum_{\mathfrak{f}\in\mathcal{F} \cap[0,\tau]} c_{k,\mathfrak{f}}\delta_f$ and absolute
convergence with respect to the total variation holds in \eqref{fact1},
because $\mathcal{F}\cap[0,\tau]$ is a finite set. If $\tau<0$, then obviously
$\mu_{k,\tau}=0$. Moreover, as  the constant $K_{\mathfrak{f}}$ in \eqref{prop_majoration_coeff_HTFc} is majorized
by an affine function of
the H\"older coefficient of $\frac{d}{dt}\mathfrak{C}_{\mathfrak{f}}$, we get from
\eqref{majoration_sol_fond_equation} that
\begin{eqnarray}
\label{eq_es_fourier_dis_sensc}
\|\mu_{k,\tau}\|=\sum_{\mathfrak{f}\in \mathcal{F} \cap[0,\tau]}\vertiii{c_{k,\mathfrak{f}}}
  \le \frac{K_1 e^{\gamma \tau}}{1+|k|^{1+ \delta}},\qquad \tau\geq0,
\end{eqnarray}
for some constant $K_1$ depending only on the $D_j$ and the $\tau_j$.
Let us now define
\begin{equation}
  \label{defmu}
\mu_k:=\sum_{\mathfrak{f}\in\mathcal{F}}c_{k,\mathfrak{f}}\delta_{\mathfrak{f}}f.
\end{equation}
Clearly, $\mu_k$ is a distribution on $\corpsR$ valued in $\corpsValeurs^{d\times d}$ and supported on $[0,+\infty)$, but it is generally not a measure. However, its restriction to
every interval bounded from above is a finite discrete measure. Hence,
we can integrate
against $\mu_k$ any $\corpsValeurs^d$-valued bounded  function which is zero for large arguments, and more generally any function that decays in norm as fast as $e^{-\gamma' t}$
for some $\gamma'>\gamma$,
by \eqref{eq_es_fourier_dis_sensc}.
Furthermore, when this function is supported on an interval of the form $[a,\tau]$ with
$a<0$ and $\tau<+\infty$, the integral against $\mu_k$ coincides with the
integral against $\mu_{k,\tau}$. 
Using this remark, we can write for $\Re(p)>\gamma$:
\begin{eqnarray}
  \label{CFG}
  G(t,p)&=&-\lim_{b\to+\infty}\int_{-\infty}^{b^+}d_\tau X(t,t-\tau)e^{-p\tau}=
            -\lim_{b\to+\infty}\sum_{k\in\mathbb{Z}}e^{i\omega kt} \int_{-\infty}^{b^+}d\mu_{k,b}(\tau)e^{-p\tau}\nonumber\\
  &=&-\sum_{k\in\mathbb{Z}}e^{i\omega kt} \lim_{b\to+\infty}\int_{-\infty}^{b^+}d\mu_{k}(\tau)e^{-p\tau}=-\sum_{k\in\mathbb{Z}}e^{i\omega kt}\int_{-\infty}^{+\infty}d\mu_k(\tau)e^{-p\tau},\nonumber  
  \end{eqnarray}
  where the second equality uses the absolute convergence in \eqref{fact1}
  and the third uses \eqref{eq_es_fourier_dis_sensc} and the fact that
  $\Re(p)>\gamma$. Considering \eqref{fourier_serie_ITF},
  the previous chain of equalities yields that
\begin{equation}
  \label{CFexp}
G_k(p)=-\int_{0^-}^{+\infty}d\mu_k(\tau)e^{-p\tau}.
\end{equation}
Next, assuming for a while that $u$ is locally bounded,
\eqref{inp-output_equation} and \eqref{fact1} imply
$$y(t) 
= \int_{-\infty}^{t^+} \sum\limits_{k \in
  \mathbb{Z}}d\mu_{k,t}(\alpha)e^{i \omega k t} u(t-\alpha)\,;$$
swapping the sum and
integral sign thanks to \eqref{eq_es_fourier_dis_sensc}
and using the above mentioned relation between $\mu_{k,t}$ and
$\mu_{k}$ yields that
$y(t) = \sum\limits_{k \in \mathbb{Z}}\int_{-\infty}^{t^+} d\mu_k(\alpha) e^{i \omega k t} u(t-\alpha)$,
and since $u(t-\alpha)=0$ when $\alpha>t$
while $\mu_k$ is supported on $[0,\infty)$ we finally get
\begin{equation}
  \label{eq:Y=GU}
  y(t) = \sum_{k \in \mathbb{Z}}A_k(t)\ \ \ \text{with}\ \ \ A_k(t):=\int_{0^-}^{+\infty} d\mu_k(\alpha) e^{i \omega k \alpha} u(t-\alpha)e^{i \omega k (t-\alpha)}\,.
\end{equation}
Since for  $\gamma'>\gamma$ we have
that
\begin{align}
  \label{estFub}
  &\!\!\!\!\int_{0^-}^{+\infty} d|\mu_k|(\alpha) \int_0^{+\infty}|u(t-\alpha)|e^{-\gamma't}dt=
    \left(\int_{0^-}^{+\infty} \!\!\!d|\mu_k|(\alpha)e^{-\gamma'\alpha}\right)\left( \int_0^{+\infty}\!\!\!|u(\tau)|e^{-\gamma'\tau}d\tau\right)
    \nonumber\\
  &  \qquad\qquad\qquad\qquad\qquad\qquad\qquad\qquad\leq\!\left(\frac{e^\gamma K_1}{1+|k|^{1+\delta}} \sum_{j=0}^{+\infty}e^{(\gamma-\gamma')j}\right)\!\!\frac{\|u\|_{L^2(\corpsR)}}{\sqrt{2\gamma'}}<+\infty 
\end{align}
%
%
%
where the  inequality uses  \eqref{eq_es_fourier_dis_sensc} and the Schwarz inequality,
we can use Fubini's theorem to compute the Laplace transform of $A_k(t)$
in the strip $\Re(p)>\gamma$. This gives us
\begin{eqnarray}
    \int_0^{+\infty}A_k(t)e^{-pt}dt&=&\int_{0^-}^{+\infty}d\mu_k(\alpha)e^{(i\omega k-p)\alpha}\int_0^{+\infty}e^{(i\omega k-p)(t-\alpha)}u(t-\alpha)dt\nonumber\\
  &=&
      \left(\int_{0^-}^{+\infty}d\mu_k(\alpha)e^{(i\omega k-p)\alpha}\right)\left(\int_0^{+\infty}e^{(i\omega k-p)\tau}u(\tau)d\tau\right)\nonumber\\
  &=&G_k(p-i\omega k)\,\hat{U}(p-i\omega k)
  \end{eqnarray}
 where the last equality uses  \eqref{CFexp}. Summing over $k$,  we find  in view of \eqref{eq:Y=GU} that
\begin{eqnarray}
\label{eq_lien_op_H_et_fourierITFc}
\hat{Y}(p)=\sum\limits_{k\in\mathbb{Z}}G_{k}(p- i \omega k )\hat{U}(p-ik\omega),\qquad\Re(p)>\gamma,
\end{eqnarray} 
where the absolute convergence of the series follows from
the estimates obtained in \eqref{estFub}; {\it i.e.},
\[
  |G_{k}(z)|\leq \frac{e^\gamma K_1K_2(\Re(z)-\gamma)}{1+|k|^{1+\delta}} \qquad\mathrm{and}\qquad|\hat{U}(z)|
  \leq \frac{\|u\|_{L^2(\corpsR)}}{\sqrt{2\Re(z)}},\qquad \Re(z)>\gamma,
  \]
  with $K_2(x):=\sum_{j=0}^{+\infty}e^{-xj}$ for $x>0$.
  Changing $p$ into $p+in\omega$ in \eqref{eq_lien_op_H_et_fourierITFc} and
  renumbering, we get 
 \begin{eqnarray}
\label{eq_lien_op_H_et_fourierITF}
\hat{Y}(p+in\omega)=\sum\limits_{m\in\mathbb{Z}}G_{n-m}(p+ i \omega m )\hat{U}(p+im\omega),\qquad n\in\mathbb{Z},\quad\Re(p)>\gamma.
\end{eqnarray} 
Combining \eqref{syst_ref1} and \eqref{eq_lien_op_H_et_fourierITF}, we see that
\begin{equation}
  \label{comparGH}
  \sum\limits_{m\in\mathbb{Z}}\left(G_{n-m}(p+ i \omega m )-R(p)^{-1}_{n,-m}\right)\hat{U}(p+im\omega)=0,
  n\in\mathbb{Z},\quad\Re(p)>\gamma,
  \end{equation}
where
$R(p)^{-1}_{i,j}$ indicates, as in \eqref{repind}, the entry with index $(i,j)$ of $R(p)^{-1}$. For each $p$ with $\Re(p)>0$, the sequence
$\{z_m:=p+im\omega,\,m\in\mathbb{N}\}$  is hyperbolically separated in the right half-plane; {\it i.e.},
\[
    \left|\frac{z_m-z_j}{z_m+\overline{z}_j}\right|\geq c >0,\qquad j,m\in\mathbb{N},\quad j\neq m.
  \]
  Indeed, we have that
  \[\left|\frac{z_m-z_j}{z_m+\overline{z}_j}\right|=\left(\frac{4(\Re(p))^2}{\omega^2|m-j|^2}+1\right)^{-1/2}\geq
    \left(\frac{4(\Re(p))^2}{\omega^2}+1\right)^{-1/2}.
  \]
  Hence (see \cite[Ch. VII, Thm 1.1]{Garnett} for an equivalent statement
  on the upper half-plane),
  $(z_m)_{m\in\mathbb{N}}$ is an interpolating sequence, meaning that to each bounded sequence
  $a_m$ in $\corpsC$ there is a bounded analytic function
  $F$ in the right half-plane with $F(z_m)=a_m$. In particular, to each
  $j\in\mathbb{N}$ and $v\in \corpsValeurs^d$, there is   a $\corpsValeurs^d$-valued
  bounded analytic function $F_{j,v}$ such that
  $F_{j,v}(z_m)=0$ for $m\neq j$ and $F_{j,v}(z_j)=v$. Multiplying $F_{j,v}$ by a function without zeros in the Hardy space $\mathcal{H}^2$
  (for instance $p\mapsto (p+1)^{-1}$), we get a function $G_{j,v}\in (\mathcal{H}^2)^d$
  which is zero at $z_m$ when $m\neq j$ and a nonzero multiple of $v$ at $z_j$. 
  Now, by the Paley-Wiener theorem, there is $u\in L^2([0,\infty)),\corpsValeurs^d)$ such that $\hat{U}=G_{j,v}$, and using this collection of
  $u$ in \eqref{comparGH} as $j$ ranges over $\mathbb{N}$ and 
  $v$ over $\corpsValeurs^d$ provides us with
  \begin{equation}
    \label{idGH}
    R(p)^{-1}_{n,-m}=G_{n-m}(p+i m \omega ),\qquad n,m\in\mathbb{N},\quad \Re(p)>\gamma.
    \end{equation}
    Next, in order to link the monodromy operator to $R(p)^{-1}$
    (this is achieved in Lemma \ref{prop_lien_mon_ITFbis} further below), we shall  use periodicity to \emph{realize} the  dynamical system \eqref{syst_lin_avec_entreebis},
    that operates in  continuous time, 
as a discrete time system in state space form; see Theorem \ref{prop_rel_mon_HTF2}.
For this, we restrict the input and the output to intervals of length $T$
by setting:
\begin{equation}
  \label{iob}
  \left.
  \begin{array}{lll}
    \widetilde{u}_k(t) &:=&u(kT+t) \\
  \widetilde{y}_k (t)&:=&y(kT+t) 
\end{array}
\right\}\mbox{ for a.e.\ $t \in [0,T]$ and $k \in \mathbb{N}$},
\end{equation}
thereby giving rise to sequences $\widetilde{u}_k$ and $\widetilde{y}_k$ of functions in
$L^2([0,T])$. We introduce another sequence of functions $\widetilde{z}_k$,
this time 
in $L^2([-\tau_N,0])$, by putting $\widetilde{z}_k:=y_{kT}$,
where the notation $y_{kT}$ was defined after \eqref{syst_lin_avec_entree0};
{\it i.e.}, $\widetilde{z}_k(\theta):=y(kT+\theta)$ for
$k\in\mathbb{N}$ and a.e.\ $\theta \in [-\tau_N,0]$. The function $\widetilde{z}_k$  will be
the ``state variable'' of our discrete system at time $k$, initialized with
$\widetilde{z}_0=y_0=0$. To realize this system,
recall the monodromy operator $U_2(T,0)$ from
\eqref{eq:stabLq} (with $q=2$) and the integration with respect to $K(t,\alpha)$ from \eqref{intvv}. These allow us to define four operators  as follows:
\begin{enumerate}
\bigskip
\item[$\bullet$]
$
\begin{array}{rl}
  \widetilde{A} :  L^2([-\tau_N,0],\corpsValeurs^d) &\longrightarrow   L^2([-\tau_N,0],\corpsValeurs^d) \\
v &\longmapsto U_2(T,0)v,
\end{array}
$
\bigskip
\item[$\bullet$]
$
\begin{array}{rl}
\widetilde{B} : L^2([0,T],\mathbb{R}^d) &\longrightarrow  L^2([-\tau_N,0],\corpsValeurs^d)\\
w &\longmapsto  \int_{0^-}^{+ \infty} d_\alpha K(T,\alpha)w(\alpha),
\end{array}
$
\bigskip
\item[$\bullet$]
$
\begin{array}{rl}
\widetilde{C} : L^2([-\tau_N,0],\corpsValeurs^d) &\longrightarrow L^2([0,T],\corpsValeurs^d) \\
v  &\longmapsto \left\{\left(U_2(t,0)v\right)(0),\ t\in [0,T]\right\},
\end{array}
$
\bigskip
\item[$\bullet$]
$
\begin{array}{rl}
\widetilde{D} : L^2([0,T],\corpsValeurs^d) &\longrightarrow   L^2([0,T],\corpsValeurs^d)\\
w &\longmapsto   \left(-\int_{0}^{+\infty} d_{\alpha}X(\cdot,\alpha)w(\alpha)\right)_{\bigl|[0,T]}.
\end{array}
$
\end{enumerate}
Observe that $\widetilde{B}$ and $\widetilde{D}$ do not depend on the values of $w$ outside $[0,T]$, and that $(U_2(t,0)v)(0)$ exists for a.e.\ $t\in\corpsR$,
hence the operators $\widetilde{A}$, $\widetilde{B}$,
$\widetilde{C}$ and $\widetilde{D}$ are  well-defined and continuous.
\begin{thm}
\label{prop_rel_mon_HTF2}
Let $u\in L^2([0,+\infty),\corpsValeurs^d)$ and $y\in L^2_{loc}([0,+\infty),\corpsValeurs^d)$ be,
respectively, the input and output  of System \eqref{syst_lin_avec_entreebis}.
For $ \widetilde{u}_k$,  $\widetilde{y}_k$ and $\widetilde{z}_k$ as in \eqref{iob},
the following recursion holds:
\begin{eqnarray}
\label{equ_syt_disc}
\left\{
\begin{array}{rl}
\widetilde{z}_{k+1} &= \widetilde{A}\widetilde{z}_{k}+\widetilde{B}\widetilde{u}_k, \\
\widetilde{y}_k&= \widetilde{C}\widetilde{z}_{k}+\widetilde{D}\widetilde{u}_k,
\end{array}
\right.\qquad \widetilde{z}_0=0,\qquad  k\in\mathbb{N}.
\end{eqnarray}
\end{thm}

\begin{pf}
Applying the variation-of-constants formula \eqref{eq:variationconstante2} to
System~\eqref{syst_lin_avec_entreebis} with $s,t$ replaced by $kT,kT+t$ yields
    \begin{equation*}
      y_{kT+t}  = U_2(kT+t,kT)y_{kT}+\int_{(kT)^-}^{+\infty}d_{\alpha}K(kT+t,\alpha)u(\alpha)\,,
    \end{equation*}
and performing the change of variable $\alpha=kT+\beta$ while using $T$-periodicity we get that
\begin{equation}
  \label{eq:505}
  y_{kT+t}= U_2(t,0)y_{kT}+\int_{0^-}^{+\infty}d_{\beta}K(t,\beta)\widetilde{u}_k(\beta)\,.
\end{equation}
When $t=T$, this gives us the first equation in   \eqref{equ_syt_disc}.
Next, evaluating each term of \eqref{eq:505} at $0$ (remember these terms
 are functions of
$\theta\in[-\tau_N,0]$ and that evaluation at a fixed point is possible for a.e.\ $t$)
while using the definition of $K$ (see \eqref{intvv}), we obtain:
\[
y(t+kT)= \left(U_2(t,0)\widetilde{z}_k\right)(0)-\int_{0^-}^{+ \infty} d_{\beta}X(t
,\beta)\widetilde{u}_k(\beta),\qquad \mathrm{a.e.}\ t\in[0,T],
\]
which is the second equation in   \eqref{equ_syt_disc}.
\qed
\end{pf}

For $a:=(a_n)_{n \in \mathbb{N}}$ a sequence in a Banach space $\mathfrak{X}$,
its $z$-\emph{transform} is the formal series:
\begin{eqnarray}
\mathcal{L} \{a\}(z) := \sum\limits_{n\in \mathbb{N}}a_n z^{-n}.
\end{eqnarray}
We let $\mathfrak{X}[[z^{-1}]]$ denote the space of such power series.
For $\mathfrak{Y}$ a Banach space and
 $O:=(O_n)_{n \in \mathbb{N}}$ a sequence of bounded
operators from $\mathfrak{X}$
to $\mathfrak{Y}$,
the $z$-transform $\mathcal{L} \{O\}(z)$ acts naturally from
$\mathfrak{X}[[z^{-1}]]$ into $\mathfrak{Y}[[z^{-1}]]$, since for each
$k\in\N$ the number of terms involved in the coefficient of $z^{-k}$ when computing
the product $ (\sum\limits_{n\in \mathbb{N}}O_n z^{-n}) (\sum\limits_{n\in \mathbb{N}}a_n z^{-n})$
is finite.
Now, if we put $\overline{y}:=(\widetilde{y}_n)_{n \in \mathbb{N}}$,
$\overline{u}:=(\widetilde{u}_n)_{n \in \mathbb{N}}$ and
$\overline{z}:=(\widetilde{z}_n)_{n \in \mathbb{N}}$, we get from
\eqref{equ_syt_disc} that
\[(zI_d-\widetilde{A})\mathcal{L}\{\overline{z}\}=
\widetilde{B}\mathcal{L}\{\overline{u}\}\qquad\mathrm{and}\qquad \mathcal{L}\{\overline{y}\}=\widetilde{C}\mathcal{L}\{\overline{z}\}+\widetilde{D}\mathcal{L}\{\overline{u}\}.\]
Note, since $\widetilde{z}_0=0$, that $(zI_d-\widetilde{A})\mathcal{L}\{\overline{z}\}$ indeed belongs to $L^2([-\tau_N,0],\corpsValeurs^d)[[z^{-1}]]$. Observing
that $(zI-\widetilde{A})^{-1}=\sum_{n\in\mathbb{N}}\widetilde{A}^{n}z^{-n-1}$
in the ring of formal Laurent series with coefficients  the operators on
$L^2([-\tau_N,0],\corpsValeurs^d)$, we may write
\begin{eqnarray}
\label{z-transform-discret-syst0bis}
\mathcal{L} \{\overline{y}\}(z)=\left[\widetilde{C}(z Id-\widetilde{A})^{-1} \widetilde{B}+\widetilde{D}\right] \,\mathcal{L} \{\overline{u}\}(z).
\end{eqnarray}
Although  \eqref{z-transform-discret-syst0bis} was derived only formally, it becomes
a valid relation between vector and operator valued analytic functions
when $|z|>\max\{\vertiii{\widetilde{A}},e^\gamma\}$ with $\gamma$ as in 
Proposition~\ref{majoration-sol-fond}, because  all series involved are then normally convergent (compare \eqref{majexpL2py}).
In another connection, writing after a change of variable the right hand side of
\eqref{inp-output_equation} as a sum of integrals on subintervals while taking into account that $X(\tau,.)$ is left continuous and using $T$-periodicity, we
get that
\begin{eqnarray}
\label{syst_discret}
  \widetilde{y}_n(t) =y(t+nT)=\sum\limits_{k=0}^{n-1}\int_{(kT)^-}^{((k+1)T)^-}
  d_{\alpha}X(t+nT,\alpha)u(\alpha)+\int_{(nT)^-}^{+\infty} d_{\alpha}X(t+nT,\alpha)u(\alpha)
  \nonumber\\
  =\sum\limits_{k=0}^{n-1}\int_{0^-}^{T^-}  d_{\beta}X(t+nT,kT+\beta)u(kT+\beta)
  +\int_{(nT)^-}^{((n+1)T)^+} d_{\alpha}X(t+nT,\alpha)u(\alpha)\nonumber\\
  =\sum\limits_{k=0}^{n-1}\int_{0^-}^{T^+}  d_{\beta}X(t+nT,kT+\beta)u(kT+\beta)
  +\int_{0^-}^{T^+} d_{\beta}X(t+nT,nT+\beta)u(nT+\beta)\nonumber\\
  =\sum\limits_{k=0}^{n-1}\int_{0^-}^{T^+}  d_{\beta}X(t+(n-k)T,\beta)\widetilde{u}_k(\beta)
  +\int_{0^-}^{T^+} d_{\beta}X(t,\beta)\widetilde{u}_n(\beta)\nonumber\\
  = \sum\limits_{k=0}^n H_{[n-k]}\widetilde{u}_{k}(t),\qquad \mathrm{a.e.}\ t\in[0,T]
\end{eqnarray}
where, for $v\in L^2([0,T],\corpsValeurs^d)$,  we have set:
\begin{eqnarray}
\label{eq_52}
H_{[k]}v(t):= \int_{0^-}^{T^+} d_{\tau}X(kT+t, \tau) v(\tau),\qquad \mathrm{a.e.} \ t\in[0,T].
\end{eqnarray}
Observe that the equality 
$\int_{(nT)^-}^{+\infty} d_{\alpha}X(t+nT,\alpha)u(\alpha)=
\int_{(nT)^-}^{((n+1)T)^+} d_{\alpha}X(t+nT,\alpha)u(\alpha)$ used in the third identity of \eqref{syst_discret} only holds for $t<T$ in general,
because when $t=T$ and $(n+1)T\in\mathcal{F}$, then
$\int_{(nT)^-}^{((n+1)T)^+} d_{\alpha}X(t+nT,\alpha)u(\alpha)$
will miss the jump that
may occur at $(n+1)T$. However, since \eqref{syst_discret} is claimed for a.e.\ $t\in[0,T]$ only, this is unimportant. Note also that changing the upper bound from $T^-$ to $T^+$ in the integrals summed over $k$ after the fourth equal sign of \eqref{syst_discret} is permitted, because the left continuity of
$X(\tau,.)$ implies that $\{T\}$ carries no mass of $d_{\beta}X(t+nT,nT+\beta)$.
If we take  $z$-transforms in $(\ref{syst_discret})$,
we obtain that
\begin{eqnarray}
\label{z-transfor-discret-systbis}
\mathcal{L} \{\overline{y}\}(z)=\mathcal{L} \{\overline{H}\}(z)\mathcal{L} \{\overline{u}\}(z)
\end{eqnarray}
where $\overline{H}:=(H_{[n]})_{n \in \mathbb{N}}$, and since for each $k$ we may pick
$\widetilde{u}_k$ arbitrarily in $L^2([0,T])$ and $\widetilde{u}_j\equiv0$ for $j\neq k$, we deduce from \eqref{z-transform-discret-syst0bis} and \eqref{z-transfor-discret-systbis} that 
\begin{eqnarray}\label{eq552}
\mathcal{L} \{\overline{H}\}(z)=\widetilde{C}(z Id-\widetilde{A})^{-1} \widetilde{B}+\widetilde{D}
\end{eqnarray}
The equality  \eqref{eq552} {\it a priori} holds
between formal power series in $z^{-1}$
with coefficients in the ring $\mathcal{B}(L^2([0,T],\corpsValeurs^d))$ of 
bounded linear operators on $L^2([0,T],\corpsValeurs^d)$. However,
if we let 
$a$ denote the spectral radius of $\widetilde{A}$ and set
$P_a:=\{z \in \corpsC:|z|>a \}$, then \eqref{eq552}  holds
in $\mathcal{B}(L^2([0,T],\corpsValeurs^d))$ when we substitute for $z$ a value in
$P_a$, for the series become normally convergent.
Of course, by analytic continuation,
the right hand-side extends
to an  operator-valued analytic function on the unbounded connected component of $\corpsC\setminus \mathrm{Spec}\,\widetilde{A}$, where $\mathrm{Spec}\,\widetilde{A}$ stands for the spectrum of
$\widetilde{A}=U_2(T,0)$, but this function may no longer be defined by a power series in $z^{-1}$.
Since we assumed that System \eqref{system_lin_formel} is exponentially stable, we know from  Proposition \ref{stability_monodromy}
that $\mathrm{Spec}\,\widetilde{A}$
is compactly included  the unit disk; that is to say: $a\in[0,1)$.
In particular,
$\mathcal{L}\{\overline{H}\}(e^{pT})$ is a well-defined
operator on
$L^2([0,T],\corpsValeurs^d)$ as soon as $\Re(p)>\log a/T$, where $\log a$
is strictly negative or $-\infty$.

For $p\in\corpsC$, consider the operator
$E_p:L^2([0,T],\corpsValeurs^d)\to L^2([0,T],\corpsValeurs^d)$ of pointwise multiplication by
the function $(t\mapsto e^{pt})\in L^\infty([0,T])$.
Clearly $E^{-1}_p=E_{-p}$, and it is easy to check that
\begin{equation}
  \label{normE}
  \vertiii{E_p}_2=\max\{1, e^{\Re(p)T}\}.
\end{equation}
When $\Re(p)>\log a/T$,
we also define the operator
\begin{eqnarray}
  \label{defL}
  \Lambda(p) &:& L^2([0,T],\corpsValeurs^d) \rightarrow L^2([0,T],\corpsValeurs^d),
  \nonumber\\
  \Lambda(p)&=&E_{-p}\circ\mathcal{L} \{\overline{H}\}(e^{pT})\circ E_p.
\end{eqnarray}
In addition, for $\Re(p)>\gamma_2$ with $\gamma_2$ as in \eqref{syst_ref1},  we construe $R(p)^{-1}:l^2(\mathbb{Z},\corpsValeurs^d)\to l^2(\mathbb{Z},\corpsValeurs^d)$ as an operator from $L^2([0,T],\corpsValeurs^d)$ into itself by identifying a function
$\phi\in L^2([0,T],\corpsValeurs^d)$ with $(\cdots,a_2,a_1,a_0,a_{-1},a_{-2},\cdots)$
where $a_k$ is the $k$\textsuperscript{th} Fourier coefficient of $\phi$. 
\begin{lem} 
\label{prop_lien_mon_ITFbis}
For $p\in\corpsC$ with $\Re(p)>\gamma_2$, it holds that
\begin{equation}
  \label{relconjFLg}
\Lambda(p)\phi=R(p)^{-1}\phi,\qquad \phi\in L^2([0,T],\corpsValeurs^d).
\end{equation}
\end{lem}
\begin{pf} As $\gamma_2>0>\log a/T$, 
  we see from \eqref{eq_52} that for any $v \in L^2([0,T],\mathbb{R})$ and $t \in [0,T]$ we have:
  \begin{equation}
    \label{calcHv}
    \mathcal{L} \{\overline{H}\}(z)v(t)= \sum\limits_{k=0}^{+\infty} H_{[k]}v(t) z^{-k} 
 = \sum\limits_{k=-\infty}^{+\infty} z^{-k} \int_{0^-}^{T^+} d_{\tau}X(kT+t,\tau)v(\tau).
\end{equation}
Performing the change of variable $\tau\rightarrow t-\tau$
in \eqref{eq_ITF} and computing as in \eqref{syst_discret}, we also obtain
\begin{eqnarray}
  \label{calcGH}
 G(t,p)&=& \int_{- \infty}^{+\infty} d_{\tau}X(t,\tau)e^{p(\tau-t)}  
 =\sum\limits_{k=-\infty}^{+\infty} \int_{0^-}^{T^+} d_{\tau}X(t,\tau-kT)e^{p(\tau-t-kT)} \nonumber\\
 &=& \sum\limits_{k=-\infty}^{+\infty} e^{-kpT} \left(e^{-pt}\int_{0^-}^{T^+} d_{\tau}X(t+kT,\tau)e^{p\tau} \right).
\end{eqnarray}
Multiplying \eqref{calcGH} on the right by some arbitrary $V\in\corpsValeurs^d$
and comparing
with \eqref{calcHv} where we set $z=e^{pT}$ and   $v(\tau)=e^{p\tau} V$ for $\tau\in[0,T]$, we get  that
\begin{eqnarray}
  \label{fconjHG}
  G(t,p)V=e^{-pt} \left[ \mathcal{L} \{\overline{H}\}(e^{pT}) 
    ( e^{p\,\cdot}V)
  \right](t),
  \qquad \mathrm{a.e.}\  t \in [0,T],
\end{eqnarray}
where $e^{p\,\cdot}V$ is another way of writing $v$ (the dot stands for a
dummy argument).
Now, assume for a while that $\phi\in C^{1,\alpha}([0,T],\corpsValeurs^d)$ with
$\alpha\in(0,1)$ and $\phi(0)=\phi(T)$. Write the Fourier expansion 
$\phi(t)=\sum\limits_{k\in\mathbb{Z}}a_k e^{i\omega kt}$ with $a_k\in\corpsValeurs^d$. By normal convergence
of the latter, we deduce from \eqref{fconjHG},
\eqref{fourier_serie_ITF} and  \eqref{prop_majoration_coeff_HTF}
that
\begin{eqnarray}
\label{eq:justif_calcul_fin1}
 \Lambda(p)\phi(t)&=&  e^{-pt}[\mathcal{L}\{\overline{H}\}(e^{pT}) e^{p \,\cdot} \sum\limits_{k \in \mathbb{Z}}a_ke^{i\omega k \,\cdot}](t)   
 = \sum\limits_{k \in \mathbb{Z}}e^{-pt}[\mathcal{L}\{\overline{H}\}(e^{pT}) e^{(p+ i\omega k) \,\cdot\,}\,a_k](t)\nonumber\\ 
 &=& \sum\limits_{k \in \mathbb{Z}}e^{i\omega kt}e^{-(p+i\omega k)t}[\mathcal{L}\{\overline{H}\}(e^{(p+i\omega k)T}) e^{(p+ i\omega k) \,\cdot\,}\,a_k](t)=\sum\limits_{k \in \mathbb{Z}} e^{i\omega k t} G(t,p+i \omega k)a_k\nonumber\\
 &=&\!\!\sum\limits_{k \in \mathbb{Z}} e^{i\omega k t} \left(\sum\limits_{n\in\mathbb{Z}}e^{in\omega t} G_n(p+i \omega k)\right)a_k=\sum\limits_{m \in \mathbb{Z}} \left(\sum\limits_{k \in \mathbb{Z}}  G_{m-k}(p+ i k \omega)a_k \right) e^{ i  m \omega t},
\end{eqnarray}
and so in view if \eqref{idGH} we get that
\begin{equation}
  \label{relconjFL}
\Lambda(p)\phi(\cdot)=\sum\limits_{m \in \mathbb{Z}} \left(\sum\limits_{k \in \mathbb{Z}}  R(p)^{-1}_{m,-k}a_k \right) e^{ i  m \omega \cdot}=R(p)^{-1}\phi(\cdot).
\end{equation}
%
Equation \eqref{relconjFL} yields the desired conclusion, except for the assumption that $\phi\in C^{1,\alpha}([0,T],\corpsValeurs^d)$ with $\phi(0)=\phi(T)$.
But since $\Lambda(p)$ and $R(p)^{-1}$ are continuous operators
$L^2([0,T])\to L^2([0,T])$ for $\Re(p)> \gamma_2$,
we conclude by a density argument that \eqref{relconjFLg} holds as stated.
\hfill$\square$
\end{pf}

To recap, we know from the discussion after \eqref{eq552}
that $\mathcal{L}\{\overline{H}\}(e^{pT})$ is an analytic operator-valued function
for $\Re(p)>\log a/T$, which is uniformly bounded in every half-plane
$\{p:\Re(p)>\beta\}$ if $\beta>\log a/T$.
Therefore $\Lambda(p)$ is, by \eqref{normE}
and \eqref{defL}, 
an  analytic operator-valued function
for $\Re(p)>\log a/T$,  uniformly bounded in every
strip $\{p:\beta_0>\Re(p)>\beta\}$ if $\beta_0>\beta>\log a/T$.
Furthermore, in view of \eqref{defHc},
$R(p)$ is an entire operator-valued function.
As Lemma \ref{prop_lien_mon_ITFbis} tells us that
$R(p)\Lambda(p)=\Lambda(p)R(p)=I_\infty$ for $\Re(p)>\gamma_2$, this identity
in fact holds on 
$\{p:\Re(p)>\log a/T\}$,
by analytic continuation.  So, if we let $\beta\in(\log a/T,0)$ and
$\beta_1=\gamma_1+1$ with $\gamma_1$ as in \eqref{VNeu},
we get that $R(p)^{-1}$ exists and is equal to $\Lambda(p)$ for
$\Re(p)>\log a/T$, and that $\vertiii{R(p)^{-1}}_2\leq C_1$
for $\beta<\Re(p)<\beta_1$ with a
constant $C_1>0$ depending on $\beta$ and $\beta_1$ but not on $p$.
In another connection, we know from \eqref{VNeu} that
$\vertiii{R(p)^{-1}}_2\leq C_2$ for
$\Re(p)>\gamma_1$, and altogether we conclude that
 

\begin{eqnarray}
\displaystyle\vertiii{R(p)^{-1}}_2\le C,\qquad \Re(p)>\beta,
\end{eqnarray}
with $C=\max\{C_1,C_2\}$. This achieves the proof of the necessity part
of Theorem~\ref{theorem_hale_generalization}, and thus of the theorem itself
in view of Section~\ref{sec:sufficiency}.

\section{Concluding remarks on neutral functional differential equations}
\label{sec:conj}
Having proven a generalization of the Henry-Hale theorem
in the periodic case,
it is natural to ask if the result from \cite{Henry1974} on the stability
of neutral functional differential systems carries over to the periodic case as well. In other words,
whether exponential stability of a periodic linear neutral differential system of the form
\eqref{eq:2} can be characterized through analyticity of its harmonic transfer function. 
It is transparent how to define the latter: simply put $\mathbfsf{H}(p):= \tilde{R}(p)^{-1}$
where
\begin{eqnarray}
\tilde{R}(p):=D_{\omega}(p)\left[I_{\infty}- \sum\limits_{j=1}^N e^{- p \tau_j} L_{D_j} \Delta_{\tau_j,\omega}\right]-\sum\limits_{k=0}^N e^{- p \tau_k} L_{B_k} \tilde{D}_{\tau_k},
\end{eqnarray}
with $\tau_0:=0$ and 
 \begin{eqnarray}
 D_{\omega}(p):=\mathrm{diag}\,\{\cdots,(p+i \omega)I_d,pI_d,(p-i \omega)I_d,\cdots\}
 \end{eqnarray}
 while the $D_j$, $B_j$ are as in \eqref{eq:2}. Then,
 with the natural definition of exponential stability
 paralleling  Definition  \ref{def:stab}, we raise the following question:

 \begin{quote}
   Assume that the $D_j(\cdot)$ and $B_k(\cdot)$ are periodic and 
differentiable with  H\"older continuous derivative for $1\leq j\leq N$ and $0\leq k\leq N$.  Is it true that
a necessary and sufficient condition for System \eqref{eq:2}
to be 
exponentially stable is the existence 
of a real 
number $\beta<0$ satisfying
\begin{enumerate}[label=\textit{(\roman*)}] 
\item \label{assumption1-conj}
$\tilde{R}(p)$ 
is invertible $\ell^2(\mathbb{Z},\corpsValeurs^d) \to \ell^2(\mathbb{Z},\corpsValeurs^d)$ for all $p$ in $\{ z \in \corpsC|\Re(z)\ge\beta\}$,
\item \label{assumption2-conj}
there is $M>0$ such that 
  $\displaystyle\vertiii{\tilde{R}(p)^{-1}}_2 \leq M$ 
  for all $p$ in $\{ z \in \corpsC|\Re(z)\ge\beta\}$ ?
\end{enumerate}
\end{quote}
We expect the answer to be positive, and it should transpire already
that necessity can be proven as
in theorem~\ref{theorem_hale_generalization}.
The difficulty with  the sufficiency part is that the 
matrix $D_{\omega}(p)$ now
prevents us from deriving straightforward   analogs of Lemmas \ref{lemma1} and \ref{lemma2}.



\end{document}